\documentclass[preprint]{article}
\usepackage[english]{babel}
\usepackage{graphics}
\usepackage{caption}
\usepackage{epsfig}
\usepackage{amssymb}
\usepackage{amsmath}
\usepackage{amsfonts}
\usepackage{float}
\usepackage{graphics}
\usepackage{graphicx}
\usepackage{multirow}
\usepackage{booktabs}
\usepackage{pstricks,pst-plot}
\usepackage{hyperref}
\usepackage{epstopdf}
\usepackage{subfigure}
\usepackage{caption}

\long\def\ignora#1{\vskip0pt}
\newcommand{\eps}{\varepsilon}

\newcommand{\diag}{\rm diag}

\newcommand{\bm}[1]{\mbox{\boldmath $#1$}}
\usepackage{listings}

\newenvironment{code1}{%
                           \mathcode`\:="603A  
                           
                           \par
                           \upshape
                           \begin{list} 
                              {} {\leftmargin = 0.0cm}
                           \item[]
                           \begin{tabbing}
                              \hspace*{.3in} \= \hspace*{.3in} \=
                              \hspace*{.3in} \= \hspace*{.3in} \=
                                                                                                                                                                                                     \hspace*{.3in} \= \hspace*{.3in} \= \kill
                          }{\end{tabbing}\end{list}}

\begin{document}



\title{A multi-class approach for ranking graph nodes: models and
experiments with incomplete data\thanks{This work was partially supported by
GNCS-INDAM,  grant ''Equazioni e funzioni di Matrici''}}


\author{Gianna M. Del Corso, Francesco Romani \thanks{Universit\`a di Pisa, Dipartimento di Informatica, Largo Pontecorvo, 3, 56127 Pisa, Italy, email: \{gianna.delcorso, francesco.romani\}@unipi.it}}
\date{April 18, 2015}
\maketitle
\begin{abstract}

After the phenomenal success of the PageRank algorithm, many re\-sear\-chers have
extended the PageRank approach to ranking graphs with richer structures
beside the simple linkage structure. In some scenarios we have to deal with
multi-parameters data where each node has additional features and there are
relationships between such features.

This paper stems from the need of a systematic approach when dealing with
multi-parameter data. We propose models and ranking algorithms which can be
used with little adjustments for a large variety of networks (bibliographic
data, patent data, twitter and social data, healthcare data). In this paper
we focus on several aspects which have not been addressed in the literature:
(1) we propose different models for ranking multi-parameters data and a class
of numerical algorithms for efficiently computing the ranking score of such
models, (2) by analyzing the stability and convergence properties of the
numerical schemes we tune a fast and stable technique for the ranking
problem, (3) we consider the issue of the robustness of our models when data
are incomplete.
%
The comparison of the rank on the incomplete data with the rank on the full
structure shows that our models compute consistent rankings whose correlation
is up to 60\% when just 10\% of the links of the attributes are maintained
suggesting the suitability of our model also when the data are incomplete.


%
%

\end{abstract}



\section{Introduction}
\label{int}

Ranking algorithms are essential tools for searching in large collections of
data and without them it would be extremely difficult to find the desired
information. Following the introduction and the success of PageRank and
similar ranking algorithms~\cite{BP98,Kl99}, researchers have extended such
techniques to a multitude of domains~\cite{BESS12,BDR,CHT09,DR09,KB06}.

In this paper we consider the setting in which the data consist of a
collection of linked items, where each item has a set of additional
attributes (features). In this setting we assume that the ranking of items
with common attributes are mutually influenced. Many important problems are
instances of this general framework. In {\em bibliographic ranking} items are
scholar papers and their citations give the linkage structure. For each paper
its associated features are its authors, the journal where it appears,
subject classification and so on. In {\em patent data} items are patents
linked by the citations to older patents. To each patent we can associate
inventors, firm, examiner, technologies, {\em etc.}. Other examples are
social or twitter graphs --- where we have information about the status, the
geographical location, the education, {\em etc.}, of users. In {\em
healthcare  data} we have patients, doctors, treatments, diseases, {\em
etc.}. With a little abuse of notation in the following we informally use the
term ``multigraph'' to denote this kind of relationships between items and
features, while other authors identify this kind of graph with as
heterogeneous information networks~\cite{SuHa12}.

In this paper we describe different models for representing the multigraph
structure of a network, and analyze different techniques for assigning
weights to features and to use these weights in the ranking process. These
weights  capture the importance that each link confers to the linked object.
We then build a fast and stable numerical method for computing the ranking
score according to our models. The proposed algorithm is obtained by
combining two non-stationary methods (BCGStab \cite{Sa03} and TFQMR
\cite{Sa03}) and a final phase of iterative refinement.

We perform many tests on two large datasets of patent data extracted from the
US patent office: the first dataset consists of all the patents granted in
the period 1976--1990 (roughly 2.5 Million patents), and the second of those
issued between 1976 and 2012 (almost 8 Million patents). The experiments aim
at understanding the role of the parameters involved in the algorithm and the
differences between the various models while comparing the results with those
returned by PageRank and the ordering induced by the citation count.

We briefly investigate also the robustness of our models when data are
incomplete and unrecoverable.
In this setting our goal is to use all the information available without
advantaging players (items or features) with more complete data respect to
those where some information is missing. We treat unknown values as zeroes,
in the sense that we do not distinguish between missing (not available) or
absent (not existent) features.
This choice is the simplest one and the one implemented in patent
repositories and in many citation databases such as {\em Scopus, Mathscinet}
where, for example, a citation is not attributed to anyone when the name of
an author has been misspelled. To evaluate the robustness of our ranking
schemas on possibly incomplete data, following the approach in related
literature~\cite{ZLZ11,GS09}, we randomly remove features from items with an
assigned probability. Our experiments show that, even removing up to half of
the features, the ranks provided by our algorithm highly correlate to the
ranks computed on the complete data. As expected, as more and more features
are removed, the ranks converge to the rank obtained using only the linkage
structure.

Finally, we tested the robustness of our models with respect to the
granularity of the features. For example if we are dealing with bibliographic
data we can  group papers into subject classes where the granularity can be
subject macro areas (Math, Computer Science, {\em etc.}) or finer
classifications (Algebra, Number Theory, Calculus, Algorithms, Data Bases,
{\em etc.}). In this context it is desirable that, when using a finer
classification, the sum of the ranks of topic A subtopics is close to A's
rank computed using the coarser classification. Experiments with the US
patent dataset show that most of our models have such desirable features.

%

The paper is organized as follows. In Section~\ref{prel} we formally
introduce the problem we are considering in the paper;  in  Section~\ref{rw}
we motivate our study and connect the techniques and the algorithm we propose
with the existing literature. In  Section~\ref{mod} we briefly present some
models discussing how extra information and features can be added to the
citation structure to improve ranking and possible weighting criteria for
such features. In our models the ranking is obtained approximating the Perron
vector of a suitable stochastic matrix.

In Section~\ref{resolution} we discuss different ways for approximating the
Perron vector showing that it can be obtained computing the solution of a
linear system. In Section~\ref{solsys} we discuss different methods for the
numerical solution of such linear system and we describe the databases used
for the experiments. In Section~\ref{numexp} we report an extensive numerical
testing to compare the different models in terms of convergence for missing
data and consistence for class aggregation. Section~\ref{conclusions}
contains the conclusion and some discussions about possible improvements of
the models.

\subsection{Preliminaries and notations on multigraphs}\label{prel}

In this paper we consider a multigraph as described by a directed graph
$G=(V, E)$ and  two mapping functions, one for the nodes $\tau: V\to {\cal
A}$ and one for the edges $\phi: E\to {\cal R}$. Each node $v \in V$  belongs
to a particular type $\tau(v)\in {\cal A}$ and each edge $e \in E$ belongs to
a particular type of relation $\phi(e)\in {\cal R}$. Functions $\phi$ and
$\tau$ are such that if $e_1$ and $e_2$ are two edges, $e_1=(v_1, v_2)$ and
$e_2=(w_1, w_2)$, with $\phi(e_1)=\phi(e_2)$, then  $\tau(v_1)=\tau(w_1)$ and
$\tau(v_2)=\tau(w_2)$. When $|{\cal A} |>1$ and ${|\cal R|}>1$  we say that
the graph is a multigraph.

\begin{figure}[htb!]
\centering%
\includegraphics[width=7cm]{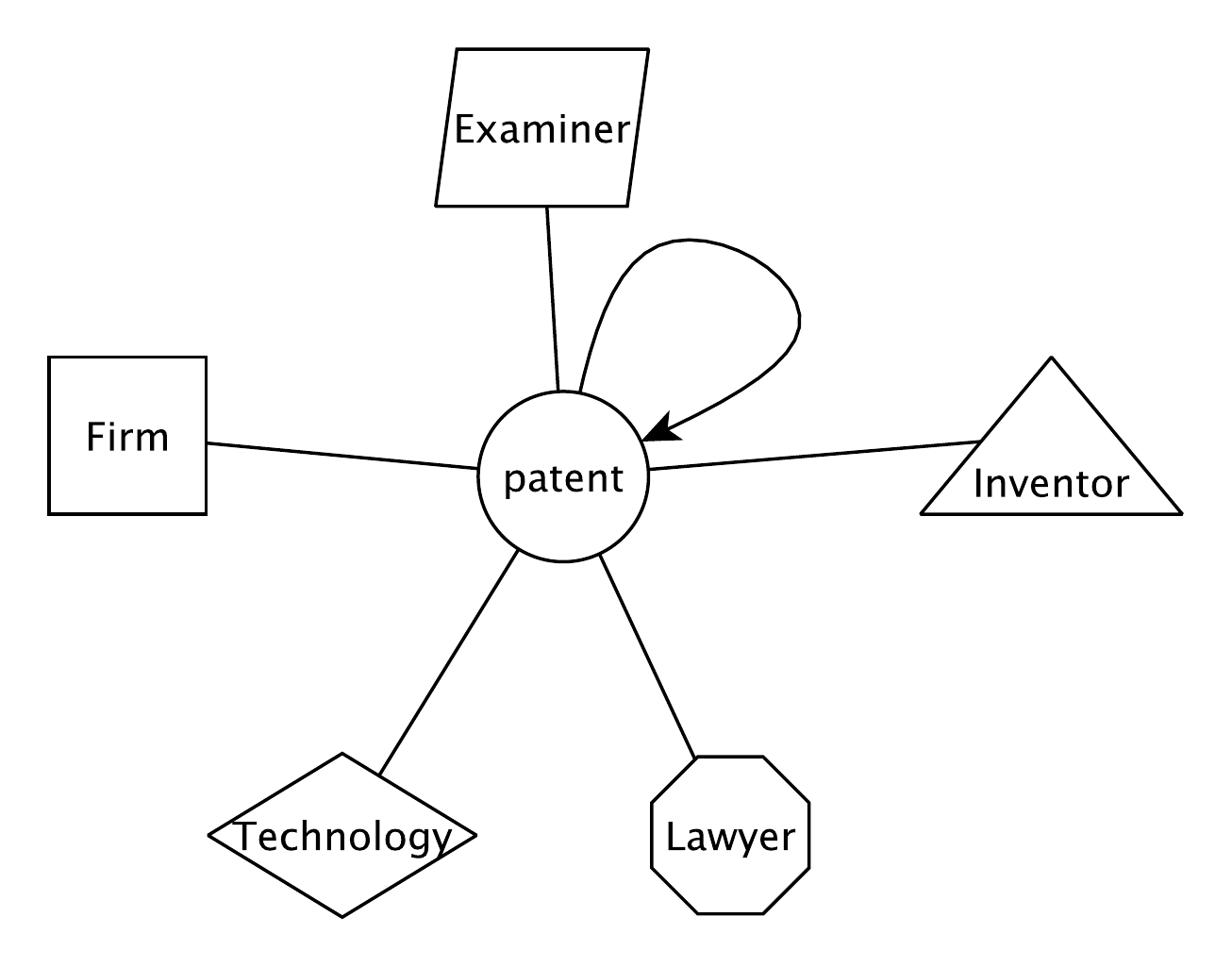}
\caption{Schema of a patent graph. Each patent has a relation with other types of nodes. } \label{heterogeneous}
\end{figure}

A typical multigraph is a patent network, where each Patent has associated  different  features in the set ${\cal A}=\left\{\right.$Patent, Technology, Firm, Examiner, Inventor and Lawyer$\left. \right\}$. The different relation types are the edges between patents and firms, patents and examiners, patents and the set of inventors, and between patents and  lawyers, beside to the edges to other cited patents: each kind of edge has a different semantic meaning, for example the connection between inventors and patents expresses intellectual property over the patent while the edge between patent and examiner represent the fact that a patent was granted by a particular examiner. In Figure~\ref{heterogeneous} it is shown the relations between the different features, and the different kinds of nodes.

To better understand multigraphs and the information contained and expressed by the different types of relations (edges) between nodes, we associate to the graph a model describing the interaction between items, features and the possible interactions between features. In this paper we define several models and compare them on the basis of a ranking function inspired by the PageRank algorithm. In particular the original network schema is enriched by including in the model other information that can be derived from the relations between nodes, such as the network of co-inventors, or all the combinations between any two couples of features, i.e. firm-technology or examiner-inventor, {\it etc}.
These enriched models allow to define a ranking function mapping objects to a real non-negative score representing  the importance of the object. The rank accounts of all the information available and not only of the citation network, and allows to rank all types of nodes on the basis of the linkage structure of the enriched graph.

\subsection{Motivations and Related work}\label{rw}

Ranking algorithms are essential when searching in large collections of data,
being either web pages, bibliographic items or even healthcare data.
Recently, many ranking algorithms have been
developed~\cite{BP98,Kl99,NZW05,Su08} which take advantage of the specific
structure of the underlined graph. Also in the area of economics it is common
practice to use ranking metrics for evaluating the performance of markets and
country economies. Recently, \cite{Monte1} has proposed a ranking algorithm
based on PageRank for patent data. Despite the whole information about
patents is available from the USPO (US Patent office) only the citation
structure has been considered in~\cite{Monte1} and the multigraph structure
of the patent graph, including also information on firms, inventors and
technologies, has not been fully exploited. Since patents are often used to
measure innovation of entrepreneurial activities~\cite{Ar62,Co02} a ranking
schema taking into account all the features of patents can be used not only
for evaluating the innovation of the patented idea or product, but also to
evaluate firms or for portfolio management. This is the primary reason we
tested the ranking algorithm presented in this paper on patent data, even if
the technique we present can be applied to any multigraph structure.


Comparing different ranking algorithms is a very difficult task since for
this problem no golden truth is available. In some cases it is possible to
take a panel of volunteers and let them manually evaluate the data, but in
most cases, either for the size of the data or for the expertise required,
this is not possible. For example, manually ranking patents requires a
remarkable knowledge of the field and such expertise is not easy to find.
Another difficulty in comparing ranking algorithms is that we can use the
same data to discover different properties: in this case a direct comparison
is not possible. For instance, if we want to evaluate scholars on the basis
of their ability to work in a team, we will design a ranking function highly
valuing the scholars with many coauthors, while if we are interested in
scientific personal strength it is natural to normalize each publication by
the number of coauthors. The resulting rankings will be completely
incommensurable.

In this paper we propose a tunable ranking algorithm where by changing
parameters we can accomplish different goals. In particular, the same
algorithm can be used on different kind of data and for different purposes.
One of the parameters is the model itself and another one is the weighting
strategy. This is the major difference with previous ranking algorithms which
are designed for specific networks and appear to be less
tunable~\cite{Su08,NZW05,SuHa12,SuYuHa09}. Together with the models we
propose and analyze some weighting strategies. To change the ranking function
one can implement other weighting strategies and incorporate them into the
algorithm.

In the following we review other approaches for ranking multigraphs and
compare them with our strategy. The problem of ranking ``multigraphs'', as
informally defined above, has been recently considered in some specific
domains. In~\cite{Su08} the multigraph is transformed into a layered graph
with a layer for each feature. The ranks of each layer are computed
independently and the final ranks are obtained with a linear combination of
the layer ranks'. We believe the independent computation for each layer does
not fully take advantage of the structure of the problem. For the specific
domain of bibliographic ranking, the PopRank algorithm powering Microsoft
Academic Search introduced by Nie et al. \cite{NZW05} is a two phase
extension of PageRank applied to typed multigraphs with different weights on
the links. In particular the formula for the PopRank score combines with
weight $\eps$ the so called ``web popularity'' which is a measure similar to
the PageRank and with weight $1-\eps$ the popularity propagation factor of
ongoing links. This factor is based on the importance of links pointing to an
object and is computed with a learning based technique which automatically
learn the popularity propagation factor for different types of links using
the partial ranking of the objects given by {\em domain experts}. This
ranking schema is very different from ours since PopRank uses an external
human contribution and is therefore problem dependent and impossible to
replicate on a different dataset.

A different approach for ranking multigraphs is the one which makes use of
multilinear algebra and tensors for representing graphs with multiple
linkages~\cite{KB06,DKK06,ADKM11}. The tensor however does not contain the
same information we use in this paper. For example, if we are dealing with
bibliographic data, our models use the full author list for each paper, while
the tensor only records the number of common authors between each pair of
papers. Hence it does not allow to obtain a score for all the features such
as authors or journals, and hence it is not possible to compare its results
with those provided by our algorithm.

Sun and al.~\cite{SuHa12,SuYuHa09} in the context of a bi-typed network (for
example a bipartite bibliographic graph with only authors and conference
venues) or star-typed networks (for example a bibliographic graph where we
have papers and all the other features such as authors, conference venues,
terms, are linked via papers)  propose a ranking schema combined with
clustering, where the clustering algorithm improves the ranking and
vice-versa. One of the ranking function proposed is similar to ours but
applies only to the simpler graphs described above with only two types of
nodes. In~\cite{ZaFe11} the authors, still in the context of bibliographic
data, proposed a model similar to one of our models, namely the {\tt Simple
Heap} model~\eqref{h1}. It mainly differs from ours for the weighting
strategy and the use of a non-static model. However we consider an enriched
structure with a complete set of relations between features. For example in
the contest of bibliographic data we enrich the graph adding weighed links
between authors, journals, and subject classification.

In previous papers from the same authors~\cite{BDR,DR09,BDR2} a model is
introduced in the contest of bibliographic data which is similar to one of
the models (the one we called {\tt Stiff model}) of this paper. In particular
in~\cite{BDR} an integrated model for ranking scientific publications
together with authors and journals was presented. In that context, particular
weighting strategies were implemented~\cite{BDR2} and an exponential decay
factor was introduced~\cite{DR09} to take into account aging of citations,
i.e. the fact that if an old  paper is not cited recently its importance
should fade over the time. In this paper we further generalize the original
ideas introducing several models and other classes making the model suitable
also for ranking other multi parameters data (patents, healthcare, social
data {\em etc.}). The new models are more adequate for example to handle
updating of the datasets which can be done at a lower cost than in the {\tt
Stiff} model. In addition, in this paper the weighting strategies are problem
independent, while in the previous papers they were designed ad hoc for
dealing with bibliographic items.


Another contribution of this paper is the investigation of adequate numerical
techniques to compute the ranking score. In particular, in
Section~\ref{resolution} we show how the computation of the ranks relies upon
the solution of a structured linear system and in Section~\ref{solsys} we
discuss and compare the different algorithms which can be used to solve that
system. Dealing with big data requires indeed particular care in the choice
of the numerical methods used in the algorithms that should be stabile and
fast. The final algorithm (Procedure {\tt SystemSolver} in
Section~\ref{solsys}) has been chosen on the basis of several tests aiming to
validate its properties of convergence and stability. A similar analysis has
not been done in the literature, and often even methods requiring matrix
manipulations~\cite{SuYuHa09} or spectral algorithms~\cite{ZaFe11} miss to
analyze this important aspect.


Another contribution of the paper, is a first analysis of the robustness of
the algorithms in the presence of missing data. Many real-world data have
missing entries and many techniques have been developed to deal with
incomplete data and to make it possible to use those dataset. A common
practice-- and the easiest to apply--  is to use only the items with complete
information discarding those with incomplete data~\cite{rev1}. This is a
rather drastic approach especially when a large portion of data is
incomplete. As an alternative, researchers have proposed to fill in a
plausible value for the missing observations. Among statisticians
distributional models for the data, such as maximum
likelihood~\cite{LR87,Sh97} and single or multiple
imputation~\cite{Sh97,Sc02}, have been developed to replace non ignorable
missing data.
The goal of this paper is not however to study the preprocessing of data for
recovery missing features. This topic would require adequate models and
techniques~\cite{LR87,rev1} to recover data and fill in the missing entries.
In this paper we are only interested in quantifying how the ranking score is
affected when some of the data are missing (completely) at
random\footnote{The data are missing completely at random (MCAR) when  the
probability that a data is missing cannot depend on any other data in the
model~\cite{Al09}. Alternative assumptions have been studied in the
literature~\cite{LR87,Sh97,Al09} such as the Missing at Random (MAR) or the
Not Missing at Random (NMAR) cases.}. To this end we assume that a missing
entry corresponds to a zero value in the linkage structure, such as is done
in large bibliographic databases such as {\em Scopus}, {\em dblp}, or even
the web when to a broken link we do not associate any link. We are aware that
replacing a missing value with a zero is not a good choice when the data do
not have homogeneous attributes~\cite{ZZJZX}, but in the case of
bibliographic data, patent data or other networks fitting into the model of
Figure~\ref{heterogeneous}, the set of the features is homogeneous. For
instance, any paper has at least an author, a publication venue, {\em etc.}
Adding and removing links at random is a common practice when evaluating
performance of ranking algorithms on large social networks to measure the
tolerance of ranking against spurious and missing
links~\cite{GS09,ZLZ09,ZLZ11}. In Section~\ref{convergence} we show that also
for our algorithm the ranking obtained with incomplete data highly correlates
with the ranking obtained with the full dataset. Of course, our analysis does
not rule out that in certain contexts an appropriate preprocessing for
recovering missing data can improve the ranking provided by our algorithm.


\section{Models}\label{mod}

In Section~\ref{oc} we present a link-based ranking for a simple citation
graph. In Section~\ref{multi-models} we enrich the graph with additional
information (features) on the nodes.

\subsection{The One-class model}\label{oc}

In this model we have a citation matrix $C$, where $c_{ij}=1$ if node $i$
links to node $j$. There are many example of such matrices for example the
web graph or the graph representing citations between scholar papers.

Following an idea similar to Google's PageRank~\cite{BP98}, we assume that
the importance $p_j$ of node $j$ is given by the importance of the nodes $i$
citing $j$, scaled by $d_i$, the outdegree of $i$. The importance given by
$i$ is thus uniformly distributed among all the cited nodes, and the
principle that the importance of a subject is neither destroyed nor created
is respected.

Here and below, we denote by $\bm e$ the vector of appropriate
length with all components equal to one. We denote by $\bm e_k$ the
$k$-th column of the identity matrix of appropriate size. The size of
vectors and matrices, if not specified, is deduced by the context.
Given a vector $\bm v=(v_i)$ of $n$ components, with the expression
$\diag(\bm v)$ we denote the $n\times n$ diagonal matrix having
diagonal entries $v_i$, $i=1,\ldots, n$.

Since nodes may have an empty set of links, the matrix $C$ can have some null
rows and in that case the corresponding outdegrees $d_i$ are zero. To avoid
divisions by zero we introduce a {\em dummy node}, numbered $n+1$, which
cites and is cited by all the existing nodes except itself. The new adjacency
matrix of size $n+1$, denoted by $\hat C$, has no null rows and is
irreducible.  The dummy node collects the importance of all the nodes and
redistributes them uniformly to all its neighbors.

The outdegrees $d_i=\sum_{j}\hat{c}_{i,j}$ define the vector $\bm d=(d_i)$,
which satisfies the equation $\bm d=\hat C \bm e$. Moreover, since $d_i\ne 0$
for all $i$, the matrix
\[
P=(p_{i,j})=\diag(\bm d)^{-1}\hat C
\]
is row-stochastic, that is, $0\le p_{ij}\le 1$, $\sum_j p_{i,j}=1$.

A similar approach is used in the PageRank model where $C$ is first
normalized by row, and then a random jump probability $\alpha$ is introduced
to make the matrix irreducible. In our model the probability to reach the
dummy node is not the same for all nodes, but varies accordingly with the
outdegree of each node.


The ranking or ``importance'' of each node is computed solving the following
equation
\begin{equation}
  \label{eq:1}
  \bm x^T=\bm x^T P, \qquad P=\diag(\hat{C}\bm e)^{-1}{\hat C}.
\end{equation}
Since the matrix $ \diag(P\bm e)^{-1}{\hat C}$ is nonnegative and
irreducible, from the Perron-Frobenius theorem~\cite{Va00} there exists a
unique vector $\bm x=(x_i)$ such that $x_i>0$, $\sum_{i}x_i=1$, which solves
\eqref{eq:1}. We call $\bm x$ the {\em Perron vector} of $P$.

This model, that we call {\tt One-class} has been introduced in~\cite{BDR}.
It has been used to rank scientific papers~\cite{BDR} and
patents~\cite{Monte1}. In~\cite{BF12}, assuming the citation matrix
triangular, this model and the PageRank model are viewed as special cases of
a family of Markov chain-based models.

\subsection{Multi-class models} \label{multi-models}

Often, beside the linkage structure we have additional information that can
be profitably used in the ranking process. For example, to evaluate a paper
we can use, besides the received citations, other information available such
as the authors or the journal where the paper has been published. We now show
that the mixing of all these ingredients (in this example authors, citations,
journals) makes it possible to compute a better ranking for papers and, at
the same time, a ranking score also for journals and authors.

The idea is to compute a ranking value for authors based on the quality of
their papers and of the journals where the papers appeared. Journals can be
evaluated as well using the information about the importance of the authors
writing for that journal and of the papers published therein. This approach
was first proposed in~\cite{BDR} and further extended in~\cite{DR09,BDR2}. We
start with the original citation matrix $C$, then we add the information on
the features of each item storing them in rectangular matrices. Examples of
features are authors and journals if the items are scholar papers; or firms,
inventors, technologies and lawyers if the items are patents. In general, we
have $f$, $f={|\cal{A}|}-1$,  rectangular binary feature matrices $F_1, \ldots, F_f$ (one for each
feature) where entry $(i, j)$ in $F_k$ is different from zero iff item $i$
has attribute $j$ for feature $k$. For patent items, for example, we
have the ``inventorship feature matrix'' storing information about the inventors
of a patent, that is, entry $(i,j)$ is nonzero if $j$ is an inventor of patent
$i$.

Given the $n_C\times n_C$ citation matrix $C$, the feature matrices $F_k$,
for $k=1, 2, \ldots , f$ where each $F_k$ has size $n_C\times n_{k}$, and
some weights $\alpha_{ij}$, we can construct a block matrix $A$ of size
$N=n_C+\sum_{k=1}^f n_k$ in different ways leading to different models. Note
that the size of $A$ is equal to the number of items plus the number of
attributes for each feature.

%

Once we have the block matrix $A$, we proceed as in the PageRank algorithm
and we obtain ranks for both {\em items} and {\em attributes}. To compute the
ranking score as in~\eqref{eq:1} we first force irreducibility in the
underlying Markov chain and then normalize the resulting matrix to get a
stochastic matrix $P$.


We now show that by varying the structure of the blocks combining the
features and the strategy for forcing irreducibility we get four different
base models. Combining these base models with different weighting strategies
we obtain a total of 15 models summarized in Table~\ref{metodi}.
\begin{description}
\item{\tt Stiff model} Each matrix $F_k$ as well as the matrix $C$ is embedded in a matrix with an additional row and column as follows
$$
\hat F_k= \left[\begin{array}{c|c} F_k & {\bm e}  \\ \hline {\bm e}^T & 0\end{array}
\right], \quad \hat C= \left[\begin{array}{c|c} C & {\bm e}  \\ \hline {\bm e}^T & 0\end{array}
\right].
$$
The matrix $\hat A$ is 
\begin{equation}\label{matA}
\hat A=\left[\begin{array}{ccccc}
\hat F_1^T\,\hat C\,\hat F_1 & \hat F_1^T\, \hat F_2       &\cdots &  \hat F_1^T\,\hat F_f& \hat F_1^T\\
\hat F_2^T\, \hat F_1        & \hat F_2^T\, \hat C\, \hat F_2&\cdots &\hat F_2^T\,\hat F_f &\hat F_2^T \\
\vdots               &\ddots               & \ddots&      \cdots           & \vdots\\
\hat F_f^T \hat F_1         &\cdots                &\ddots               &\hat F_f^T\,\hat C\,\hat F_f & \hat F_f^T\\
\hat F_1      	      &\hat F_2		    &\cdots &\hat F_f		&\hat C
\end{array}
\right].	
\end{equation}
The matrix $\hat A$ is the adjacency matrix of a more complex multigraph respect to the one described by the schema in Figure~\ref{heterogeneous}. In fact all the possible relations between any pair of features is accounted for, meaning that the graph is complete and we have $ (f+1)^2$ types of edges.
The diagonal blocks are of the form $\hat F_k^T\hat C\hat F_k$ and
contain the co-citations between features. For example, if $F_k$ is the
authorship matrix each entry of $\hat F_k^T\hat C\hat F_k$ accounts for
the number of citations between any two authors. For off-diagonal blocks
of type $\hat F_k^T \hat F_h$, for example when $F_h$ is the
paper-journal matrix\footnote{In the paper-journal matrix an entry
$(i,j)$ is nonzero if the paper $i$ was published on journal $j$.},
each entry accounts for how many papers an author has published on a
given journal.

For the construction of the stochastic and irreducible matrix $P$ we
proceed as follows. We normalize by row each block of matrix $\hat A$,
obtaining the stochastic and irreducible matrices $P_{i,j}$, for $i, j=1,
2, \ldots, f, f+1$, where $P_{f+1, f+1}$ corresponds to the row
normalization of $\hat C$. Then, given a row stochastic matrix of weights
$\Gamma=(\gamma_{ij})$ with $i, j=1, 2, \ldots, f, f+1$, we build matrix
$P$ as follows
\begin{equation}\label{matP}
P=\left[\begin{array}{cccc}
\gamma_{1,1}\, P_{1,1} & \gamma_{1,2}\, P_{1,2}       &\cdots &  \gamma_{1,f+1}\, P_{1,f+1}\\
\gamma_{2,1}\, P_{2,1}      & \gamma_{2,2}\, P_{2,2}&\cdots & \gamma_{2,f+1}\, P_{2,f+1} \\
\vdots               &\vdots                        & \vdots\\
\gamma_{f,1}\, P_{f,1}         &\cdots                           &\gamma_{f,f}\, P_{f,f} & \gamma_{f,f+1}\, P_{f,f+1}\\
\gamma_{f+1,1}\, P_{f+1,1}      	      		    &\cdots & \gamma_{f+1,f}\, P_{f+1,f}		& \gamma_{f+1,f+1}\, P_{f+1,f+1}
\end{array}
\right].	
\end{equation}
We called this model {\tt Stiff} because it lacks flexibility. In fact,
if we add an attribute to a feature, we need to recompute not only the
corresponding $F_k$ and the matrices involving $F_k$ in~(\ref{matA}), but
also renormalize each of the changed blocks. This approach was followed
in~\cite{DR09,BDR2} for ranking papers, authors and
journals\footnote{In~\cite{DR09,BDR2} each block was normalized in a
particular way because row normalization was not always well suited for
that particular problem.}. In~\cite{DR09} some discussion about possible
choices of the weights $\gamma_{ij}$ are reported. Note that since  the
matrix of the weights $\Gamma$ is stochastic and also the blocks $P_{ij}$
are stochastic, matrix $P$ describes a coupled Markov chain.

\item{\tt Static model} This model differs from the previous because
    instead of adding a row and a column to each of the feature matrices,
    we add a dummy item to the whole matrix, and then weight each block
    with suitable parameters $\alpha_{ij}$. We obtain the matrix
\begin{equation} \label{static}\hat A=\left[\begin{array}{c|c}
\begin{array}{ccccc}
\alpha_{1,1}\,F_1^T\, C\, F_1 & \alpha_{1,2}\,F_1^T\, F_2       &\cdots & \alpha_{1,f}\, F_1^T\,F_f& \alpha_{1,f+1}\,F_1^T\\
\alpha_{2,1}\,F_2^T\, F_1        & \alpha_{2,2}\,F_1^T\, C\, F_1&\cdots & \alpha_{2,f}\,F_2^T\,F_f &\alpha_{2,f+11}\,F_2^T \\
\vdots               &\ddots               & \ddots&      \cdots           & \vdots\\
\alpha_{f,1}\,F_f^T F_1         &\cdots                &\ddots               &\alpha_{f,f}\,F_f^TCF_f & \alpha_{f,f+1}\,F_f^T\\
\alpha_{f+1,1}\,F_1      	      &\alpha_{f+1,2}\,F_2		    &\cdots & \alpha_{f+1,f}\,F_f		&\alpha_{f+1,f+1}\,\tilde C
\end{array} & {\bm e}\\ \hline
{\bm e}^T & 0
\end{array}
\right].	
\end{equation}
We then  normalize by row to get the stochastic irreducible matrix
$P=\diag(\hat A\bm {e})^{-1} \, \hat A$. For this and the remaining
models proposed in this section, whenever we add a new attribute to an
existing feature we have to change only the matrix of the feature
involved. Indeed, we do not
need to build the matrix $\hat A$ explicitly but all the computation can
be done using only the matrices $F_k$ and $C$.
\end{description}

The next two models are designed for dealing with problems where the feature
data is incomplete. For example in a bibliographic database where we only know the
first author of each paper. In this case, we cannot expect to compute an
accurate rank for authors, but still we would like to use the available
author information to better rank papers. The structure of the blocks is now
homogeneous  among off-diagonal and diagonal blocks so that we can ideally
consider all the features heaped in just a matrix $F$ containing all the
information on the different attributes. Matrix $F$ has size $n_C\times n_t$,
$n_t$ being the total number of attributes, for example the sum of distinct
authors, journals, {\em etc.} available. Since the Heap model can be used
also with complete data we describe the model keeping the features distinct,
knowing that the features can be squeezed in a unique matrix when the
features classes are scarcely populated.

\begin{description}

\item {\tt Heap model} The {\tt Heap} model differs from the {\tt Static}
    model in the off-diagonal blocks.  Blocks $F_k^T F_h$ are replaced by
    $F_k^TCF_h$. In the previous example where $F_k$ was the
    paper-journal matrix and $F_h$ is the paper-author matrix, the entry
    $(i,j)$ of  $F_k^TF_h$ is the number of papers author $j$ has
    published on journal $i$, while the $(i,j)$ entry of $F_k^TCF_h$ is
    the number of citations from papers written by author $j$ to all
    papers published in journal $i$.

Assigning to each block a weight $\alpha_{i,j}$, we get  the matrix $\hat A$
$${\hspace{-0.5cm}}
\hat A=\left[\begin{array}{c|c}
\begin{array}{cccccc}
\alpha_{1,1}\,F_1^T\, C\, F_1 & \alpha_{1,2}\,F_1^T \,C\, F_2       &\cdots & \alpha_{1,f}\, F_1^T\,C\,F_f& \alpha_{1,f+1}\,F_1^T\\
\alpha_{2,1}\,F_2^T\,C \,F_1        & \alpha_{2,2}\,F_1^T\, C\, F_1&\cdots & \alpha_{2,f}\,F_2^T \,C\,F_f &\alpha_{2,f+11}\,F_2^T \\
\vdots               &\ddots               & \ddots&      \cdots           & \vdots\\
\alpha_{f,1}\,F_f^T\,C\, F_1      &\cdots                &\ddots               &\alpha_{f,f}\,F_f^T\,C\,F_f & \alpha_{f,f+1}\,F_f^T\\
\alpha_{f+1,1}\,F_1      	      &\alpha_{f+1,2}\,F_2		    &\cdots & \alpha_{f+1,f}\,F_f		&\alpha_{f+1,f+1}\,\tilde C
\end{array} & {\bm e}\\ \hline
{\bm e}^T & 0
\end{array} \right].
$$
To get the stochastic matrix $P$ we just normalize $\hat A$ by row.


\item {\tt Simple Heap model} In this model we assume that there is no
    interaction between features so that cross-citations do not influence
    the rank.
\begin{equation}\label{h1} \hat A=\left[\begin{array}{c|c}
\begin{array}{cc}
O & \alpha_{1,2} \,F^T \\
\alpha_{2,1}\, F & \alpha_{2,2} \tilde C \end{array} & \bm{e}\\
\hline {\bm e}^T
\end{array} \right],
\end{equation}
where $F$ is a matrix containing all the relations between items and attributes, {\em i.e.} $F=[F1, F2, \ldots, F_f]$.
As already observed this model uses a simplified setting to deal with the case where we have incomplete data.
%
\end{description}

that as $\alpha_{f+1, f+1}\to 1$  and $\alpha_{i,j}\to 0$, for all the other values of $i, j$, the rank obtained with these models converges to the
rank obtained with the one-class model. This is however guaranteed because in
all the models,  for the limit value of $\alpha_{i,j}$, $\hat A$ collapses to
a matrix of the form
$$
\left[\begin{array}{c|c}
\begin{array}{cc}
O & O \\
O &\tilde C \end{array} & \bm{e}\\
\hline {\bm e}^T
\end{array} \right],
$$
and the rank of the items is the same (up to a scaling factor) of the one
obtained with the one-class model, while all the features will get an uniform
score.


\subsection{Weighting strategies}\label{ws}

Weighting strategies play  an important role in the tuning of the algorithm,
since by varying them we can change the relative importance of features vs
citations and consequently change the final ranking. We propose five
different weighting strategies for our models, but not all strategies can be
applied to each model, and for different models two weighting schemes may
coincide after normalization of the matrix $\hat A$.

The simplest strategy is the {\tt Uniform} ({\tt U}) one, which corresponds
to choosing $\alpha_{i,j}=1$ for each $i, j=1, \ldots, f+1$. By adopting this
weighting schema the contribution of each class (feature or citation) is
valued in the same way, independently of its size. This approach appears
adequate only when the sizes of each class are of the same order of
magnitude, otherwise we are giving a bigger role in the determination of the
ranking to scarcely populated classes.

For this reason, we also consider schemes that keep track of the size of
each class. We have different choices.

\begin{description}

\item{\tt Dimension-based (D)} We set $\alpha_{i,j}=n_j/n_C$, and
    $\alpha_{i,f+1}=1$. In this way we guarantee that the average value
    of the features are the same~\cite{DR09}, and we do not advantage
    more populated classes respect to those less populated. The weights
    are the same for each block of columns.

\item{\tt Double-Dimension-based (DD)} We have a symmetric weight matrix,
    setting $\alpha_{i,j}=\alpha_i\alpha_j$, where $\alpha_i=n_i/n_C$ is
    the normalized size of the $i$-th feature. In the case the citation
    matrix is  much larger respect to the size of $F_i$, this scheme
    gives more  importance to citations than to features.

\item{\tt Heap (H)} We set  $\alpha=(\sum_{k=1}^f n_k)/n_C$ for the first
    $f$ blocks of columns, that is $\alpha_{i, j}=\alpha$ for $i=1,
    \ldots, f+1$ and $j=1, \ldots, f$  for the blocks in the last  column
    we get $\alpha_{j, f+1}=1$, for $j=1, \ldots, f+1$. This weighting
    strategy is particularly suited for the {\tt Heap} or {\tt Simple
    Heap} model.

\item{\tt Double-Heap (HH)} In this case the weights are not the same
    along the blocks of columns but defining $\alpha=(\sum_{k=1}^f
    n_k)/n_C$, we have $\alpha_{i,j}=\alpha^2$ for $i, j=1, \ldots, f$,
    and the weights of blocks in the last column are  $\alpha_{j,
    f+1}=\alpha$, and in the last row $\alpha_{f+1, j}=\alpha $. Moreover
    $\alpha_{f+1, f+1}=1$. Also this scheme is particularly suited for
    the {\tt Heap} or {\tt Simple Heap} model since they have the same
    value in the upper left blocks. Assuming $\alpha<1$ we are giving
    again more importance to citations when determining the ranking
    scores of the other nodes.

\end{description}

While it is always possible to apply an Uniform weight to each base model, it
doesn't make sense to apply some of the weighting strategies to the {\tt
Stiff} or {\tt Static} model. In fact the {\tt H} or the  {\tt HH} weighting
techniques make sense only when the structure of diagonal and off diagonal
blocks is the same as  in the case of the {\tt Heap} or {\tt Simple-Heap} model.
Using the {\tt H}  or the {\tt HH} weighting
techniques in combination with the {\tt Heap}  model we can rewrite the matrix $\hat A$ in a more compact form collecting all the features in a unique matrix $F$. We get
$$
\hat A=\left[\begin{array}{c|c}
\begin{array}{cc}
\alpha_{1,1}\, F^T\, C\, T & \alpha_{1,2} \,F^T \\
\alpha_{2,1}\, F & \alpha_{2,2} \tilde C \end{array} & \bm{e}\\
\hline {\bm e}^T&0
\end{array} \right].
$$
 In
Table~\ref{metodi} we summarize the fifteen full models obtained combining
the four basic models, with the five weighting schemas.
\begin{table}
\begin{center}
{
  \begin{tabular}{ ||l||r|r|r|r|r|| }
    \hline
    {\tt models}$\backslash ${\tt weigths} & {\tt U} & {\tt D} & {\tt DD} & {\tt H}& {\tt HH}
    \\ \hline \hline
    {\tt Stiff} & Stiff-U&Stiff-D &-& -& - \\
 {\tt Static} & StaticU&Static-D &Static-DD &- & - \\
   {\tt Heap} & Heap-U& Heap-D&Heap-DD &Heap-H & Heap-HH \\
   {\tt Simple-Heap} &SHeap-U & SHeap-D&SHeap-DD & SHeap-H&  SHeap-HH\\
       \hline
  \end{tabular}
  \caption{{The 15 models obtained combining the basic models with the
  different weighting strategies.}}\label{metodi}}
\end{center}
\end{table}

\section{Computation of the Perron vector}\label{resolution}

In all our models to compute the rank we have to solve an eigenvector problem
involving a stochastic irreducible matrix. More precisely, we have to find
the left Perron vector ${\bm x}$ such that ${\bm x}^T={\bm x}^T \, P$, with
$P$ stochastic. We now show that the Perron vector can be computed as the
solution of a linear system involving a matrix $\widehat M$, where
$$
\widehat M=\left\{ \begin{array}{ll}
\hat C &\mbox{ if } f=1\\
\hat A & \mbox{ if } f\ge 2.
\end{array}
\right.
$$
separating the last row and column of $\widehat M$ we have
$$
\widehat M=\left[\begin{array}{cc}
M & \bm{u}\\
\bm{v}^T & 0\end{array}
\right],
$$
where $M$ has size $N\times N$ and $\bm{u}, \bm{v}$ are suitable $N$-vectors
(for the Stiff models $\bm{u}, \bm{v}$ are the last column and row of $P$
in~(\ref{matP}), while for all other models $\bm{u}=\bm{v}=\bm{e}$). The
matrix $P$ is obtained normalizing by row  $\widehat M$, that is
$P=\diag(\widehat M\, \bm{e})^{-1}\, \widehat M$. Let
$$
D=\diag(\widehat M\, \bm{e})^{-1}=
   \left[\begin{array}{cc} D(\bm{u}) & \\
  & 1/(\bm{v}^T {\bm e})\end{array}\right],
$$
where $D(\bm{u)}=\diag(M\,\bm{e}+\bm{u})^{-1}$. Setting $\bm{x}^T=(\bar{\bm
x}^T, x_{n+1})$, where $\bar{\bm x}^T$ is an $n$-vector, the equation ${\bm
x}^T={\bm x}^T \, P$ can be rewritten as
\begin{equation}\label{sistema}
\left\{\begin{array}{l}
\bar{\bm x}^T= \bar{\bm x}^T\,D({\bm u}) \, M+\displaystyle{\frac{x_{n+1}}{\bm{v}^T\,{\bm e}}} \, {\bm v}^T\\
\\
x_{n+1}=\bar{\bm x}^T D(\bm{u})\, \bm{u}.
\end{array}\right.
\end{equation}
Since we are interested in the direction of the Perron vector and not in its
norm, we can chose $x_{n+1}={\bm{v}^T\,{\bm e}}$, obtaining $\bar{\bm x}^T=
\bar{\bm x}^T \,D({\bm u}) \, M+ {\bm v}^T$. The vector $\bar{\bm x}$ is then
the solution of the linear system
\begin{equation}\label{ls}
\left(I-M^T D(\bm{u})\right){\bar {\bm{x}}}={\bm{v}},
\end{equation}
or can be computed by the iterative method
\begin{equation}\label{mi}
{\bar {\bm x}^{T\,{(i+1)}}}= {\bar {\bm x}^{T\,{(i)}}} D(\bm{u}) \, M +\bm{v}^T.
\end{equation}
Note that for the {\tt Stiff} models it is $D=I$.

It is important to observe that in the proposed models we can simply work
with the matrices $F_j$ without explicitly normalize and store the complete
matrix $\widehat M$.  For example, for the {\tt Static} model
in~\eqref{static} the $i$-th block, $i=1, \ldots, f$ of the vector
$\widehat{M}{\bm e}=M{\bf e}+{\bf u}$, used for constructing matrix $D$,  can
be computed as follows
$$
{\bf z}_i=\sum_{j\neq i}\alpha_{i,{j}}\, F_i^TF_j {\bf e}_{j}+\alpha_{i,i}F_i^TCF_i{\bf e}_i+{\bf{e}}_i, \quad i=1, \ldots, f
$$
and
$$
{\bf z}_{f+1}=\sum_{j=1}^{f} \alpha_{f+1, j} F_j {\bf  e}_j+\alpha_{f+1, f+1}\tilde C\, {\bf e}_{f+1}+ {\bf e}_{f+1}.
$$
The cost for computing ${\bf z}_i$ is linear in the number of non zeros
(denoted as {\it nnz}) of all the matrices $F_i$ and $C$, that is $O(\sum_i
nnz(F_i)+nnz(C))$ since the matrices $F_i$ are stored in a sparse format and
the cost of  multiplying  a sparse matrix by a vector is equal to the number
of non zeros in the matrix.

Letting ${\bm w}_i$ denote the vectors of length $n_i$, whose entries are the
reciprocal of the entries of the vectors ${\bm z}_i$, and  noticing that the
$i$-th diagonal block of matrix $D({\bf u})$ contains the entries of ${\bm
w}_i$, an iteration of~\eqref{sistema} becomes
$$
\small{
{\bar {\bm x}_j^{T\,{(k+1)}}}= \left\{ \begin{array}{l}
\displaystyle{\sum_{i\neq j}\alpha_{i,j} ({\bar {\bm x}_i^{T\,{(k)}}}*{\bf w}_{i}) F_i^TF_j+
({\bar {\bm x}_j^{T\,{(k)}}}*{\bf w}_{j}) F_j^T\tilde C F_j +\bm{e}_j }\\
\hfill  j\le f\\
 \\
\displaystyle{\sum_{i=1}^f\alpha_{f+1, i} ({\bar {\bm x}_{i}^{T\,{(k)}}}*{\bf w}_{j})F_i^T+ \alpha_{f+1, f+1} ({\bar {\bm x}_{f+1}^{T\,{(k)}}}*{\bf w}_{f+1}) \tilde C}+{\bm e}_{f+1} \\
\hfill   j=f+1,
\end{array}
\right.
}
$$
where  $*$  denotes the component-wise (Hadamard) product between vectors.
The component-wise products can be computed in $O(\sum_i n_i+n_C)$
multiplications,  and the total cost of computing the new vector ${\bar {\bm
x}^{T\,{(k+1)}}}$ is proportional to the number of non zeros in the matrix $[F_1,
F_2, \ldots F_f, C]+\sum_i n_i+n_C$. We can proceed analogously on the other
models.

\section{Solution of the Linear System with non-stationary methods}\label{solsys}

Once the problem of the computation of the Perron vector is reformulated as
the solution of the linear system~\eqref{ls}  we can employ the iterative
method described in~\eqref{mi} or stationary methods such as Jacobi or
Gauss-Seidel iterations, or the more promising Krylov methods. In fact, also
non-stationary methods need only the computation of matrix-vectors products
and are in general more effective than stationary ones
(see~\cite{DGR07,GZB04} for a comparison between stationary and
non-stationary methods on similar problems). Recall that to compute the
product $M{\bm x}$ then, we do not explicitly form and store the matrices $M$
and $D({\bm u})$ but we store in sparse form only the matrices of the
features $F_i$ and the citation matrix $C$.

We implemented different Krylov methods, and among them we chose the three
more performing: BCGStab, CGS, TFQMR
(see~\cite{Sa03} for the details on these methods).

To refine the final result we add a few steps of the iterations~\eqref{mi} in
accordance with the {\tt Iterative Refinement} algorithm described below. In
particular we perform some additional iterative step  until either the
distance of two successive iterations is less than {\tt tol} or we are stuck
and the vector is not changing anymore.

{{
\framebox{\parbox{4.0cm}{
\begin{code1} \label{algo1}
\noindent{\bf Procedure} {\tt Iterative Refinement }\\
{\bf Input}: ${\bf x^{(i-1)}}, {\bf x^{(i)}}, {\bf x^{(i+1)}}$, {\tt  tol}\\
\ {\bf while} $\|{\bf x}^{(i)}- {\bf x}^{(i-1)}\|<\mbox{{\tt tol}}$ {\bf or}  $\left| \|{\bf x}^{(i)}- {\bf x}^{(i-1)}\|- \|{\bf x}^{(i+1)}- {\bf x}^{(i)}\| \right|<\mbox{{\tt tol}}$ \\
\  {\bf do}
 a step of the iterative method~\eqref{mi}, $i=i+1$\\
\ {\bf endwhile}
\end{code1}}}}}

\subsection{Models Validation: Stability and Convergence}\label{stab}

To test the methods for the solution of~\eqref{ls} we constructed two
datasets with real data extracted from the $US$ patent office and we used five
features: Firms, Inventors, Technologies, Lawyers and Examiners. In
particular, we denote by $F_1$ the patent-technology matrix where entry
$(i,j)$ is one if  patent $i$ uses technology $j$;  by $F_2$ the patent-firm
matrix, recording  the firm owning the patent, by $F_3$ the patent-inventors
matrix which maps patents to inventors, by $F_4$ the patent-lawyers where
each patent is matched to the lawyers applying for the patent, and by $F_5$
the matrix where at each patent is associated  the examiners from the US
Patent Office who approved the patent. The matrix $C$ contains the citations
between patents and is almost triangular since each patent can be based only
on patents from the past.
\begin{description}
\item{DS1:} Consists of $n_C=2\, 474\, 786$ US patents from 1976-1990. Of
    these patents we have additional information that can be grouped into
    5 major features,  namely $n_1=472$ {\tt Tech\-no\-lo\-gies}, $n_2=165\,
    662$ {\tt Firms}, $n_3=965\,878$ {\tt Inventors}, $n_4=25\,341$ {\tt
    Lawyers} and $n_5=12\,817$ {\tt Examiners}, giving rise to a matrix
    $\hat A$ of size $n_C+\sum_{i=1}^5 n_i$ which is approximately of
    3.7 millions.
\begin{table}[t!]
\begin{center}
  \begin{tabular}{ ||l||r|r||r|r||r|r|| }
    \hline
    {\tt models} &{\tt BCGstab} &  & {\tt CGS} & & {\tt TFQMR} & \\ \hline
                  &it & $\log_{10}(res)$                &it & $\log_{10}(res)$ &it & $\log_{10}(res)$\\ \hline\hline
Stiff-U&	18&	-10.49	&100&-7.71&	21&	-3.90 \\ \hline
Stiff-D&	23&	-11.77	&100		&-11.20	&19	&-4.75\\ \hline
Static-U&	35&	-9.03	         &100  		&-6.25	&40	&-7.83\\ \hline
Static-D&	39&	-11.13	&100		&-7.22	&37	&-9.60\\ \hline
Static-DD&	35&	-12.33	&100		&-12.20	&30	&-11.99\\ \hline
Heap-U&	32&	-9.86		&100		&-7.44	&36	&-8.59\\ \hline
Heap-D&	36&	-11.26	&100		&-7.73	&38	&-9.74\\ \hline
Heap-DD&	41&	-11.48	&100		&-9.46	&33	&-11.51\\ \hline
Heap-H&	36&	-10.83	&100	&-6.46	&30	&-7.72\\ \hline
Heap-HH&	24&	-9.85		&100	&-7.14	&27	&-8.38\\ \hline
SHeap-U&	32&	-11.56	&100	&-8.00	&29	&-9.91\\ \hline
SHeap-D&	32&	-11.72	&100	&-8.51	&28	&-9.85\\ \hline
SHeap-DD&	37&	-11.43	&100	&-11.83	&28	&-11.97\\ \hline
SHeap-H&	28&	-10.56	&100	&-8.33	&25	&-9.98\\ \hline
SHeap-HH&	29&	-11.34	&100	&-6.75	&24	&-10.05\\ \hline
  \end{tabular}
  \caption{{Performance comparison between
  three Krylov methods on the 15 models on a problem of size 3.7 million.}}\label{tab1}
\end{center}
\end{table}

\item{DS2:} Consists of $7\, 984\, 635$ US patents from 1976-2012. The
    size of the five features are as follows $475$ {\tt Technologies},
    $633\,551$ {\tt Firms}, $4\,088\,585$ {\tt Inventors}, 120\,668 {\tt
    Lawyers} and 64\,088 {\tt Examiners}, giving rise to a matrix $\hat A$ of size
    approximately of 13 millions.
\end{description}
{\begin{table}[t!]
\begin{center}
  \begin{tabular}{ ||l||r|r||r|r||r|r| }
    \hline
     {\tt models}	&{\tt BCGstab}&	  	&{\tt TFQMR}	   & \\ \hline
               &it & $\log_{10}(res)$                &it & $\log_{10}(res)$  \\ \hline\hline
Stiff-U&	14&	 -10.57&	 19&	 -3.83  \\ \hline
Stiff-D&	21&	 -11.30&	 25&	 -3.71\\ \hline
Static-U&	 37&	 -6.77&	 52&	 -3.09\\ \hline
Static-D&	52&	 -10.53&	 53&	 -7.89\\ \hline
Static-DD&	39&	 -11.38&	 43&	 -8.70\\ \hline
Heap-U&	 36&	 -8.87&	 47&	 -7.58\\ \hline
Heap-D&	45&	 -6.47&	 51&	 -6.11\\ \hline
Heap-DD&	41&	 -9.40&	 41&	 -6.56\\ \hline
Heap-H&	 40&	 -9.63&	 43&	 -7.31\\ \hline
Heap-HH&	35&	 -9.49&	 43&	 -7.52\\ \hline
SHeap-U&	40&	 -9.78&	 38&	 -7.48\\ \hline
SHeap-D&	38&	 -10.36&	 36&	 -8.10\\ \hline
SHeap-DD&	36&	-11.75&	 34&	 -9.79\\ \hline
SHeap-H&	31&	 -7.90&	 35&	 -4.54\\ \hline
SHeap-HH&	35&	-10.63&	 35&	 -5.91\\ \hline
\end{tabular}
  \caption{{Performance comparison between two Krylov methods
  applied to the 15 models of Table~\protect{\ref{metodi}} on a problem
  of size 13 million.}} \label{tab2}
\end{center}
\end{table}
}

The feature matrices and the citation matrix $C$ are used to obtain ranks both for patents and
features, i.e. Technologies, Firms, Inventors Lawyers and Examiners with the techniques
described in this section.

When using iterative solvers we have always to address the question of
numerical stability. The three proposed methods, {\tt BCGStab}, {\tt CGS} and
{\tt  TFQMR} have been tested on the two datasets with an error goal of
$10^{-11}$ and with maximum number of iterations equal to 100. For the
refinement steps of the power method we set ${\tt tol}=10^{-13}$.
Applying to dataset DS1 the three methods to all the models we obtain the
results summarized in Table~\ref{tab1}, where instead of the actual residuals
we report only their  base 10 logarithm.

It is evident that {\tt CGS} is inadequate to cope with this kind of problems
since after 100 iterations we have still a high residual norm. Moreover {\tt
BCGstab} is better then {\tt TFQMR} since it achieves almost always a lower
residual norm. For these reasons we restrict our analysis to {\tt BCGStab}
and {\tt TFQMR} comparing them on the dataset of size 13M. We obtain the
results reported in Table~\ref{tab2}.

We note that {\tt BCGstab} is clearly better than {\tt TFQMR}, but sometimes
fails to reach an acceptable accuracy. Hence a three step algorithm,
described in Procedure {\tt SystemSolver}  has been devised.

\medskip
{{
\framebox{\parbox{8.0cm}{
\begin{code1} \label{algo2}
{\bf Procedure} {\tt SystemSolver }\\
{\bf Input}: Initial guess ${\bf x^{(0)}}$, {\tt ErrorGoal}, {\tt maxiter}, {\tt  tol}\\
\ {\bf Apply} {\tt BCGStab} with error goal={\tt ErrorGoal} and maximum iterations={\tt maxiter}  \\
\ {\bf if} {\tt res} $>$ {\tt ErrorGoal} \\
\quad   {\bf Apply} {\tt TFQMR} with error goal={\tt ErrorGoal} and maximum iterations={\tt maxiter}  \\
\ {\bf endif} \\
\  {\bf Apply} {\tt Iterative Refinement} with tolerance {\tt tol}
\end{code1}}}}}
\medskip

Applying this procedure, with {\tt ErrorGoal}=$10^{-10}$, {\tt maxiter}=100
and {\tt tol}=$10^{-13}$, on  both the datasets we get the results displayed
in Table~\ref{tab3}.

From Table~\ref{tab3} we observe that the models which are more stable for
the two datasets considered are the {\tt Stiff-D}, {\tt Static-DD}, and among
the Heap-like models, we have good performance of {\tt Heap-DD}, {\tt
SHeap-DD}.
 \begin{table}[t!]
\begin{center}
{\begin{tabular}{||l||r|r||r|r|| }\hline
{\tt models} &DS1&size=3.7M&	DS2 &	size =13M \\ \hline
&{\tt time(sec.)} &$\log_{10}(res)$	&{\tt time(sec.)}	&$\log_{10}(res)$\\ \hline \hline					
{Stiff-U}&	237&	-12.095&		1078	&-11.952\\ \hline
{Stiff-D}&	179&	-12.762	&	1422&	-12.1884\\ \hline	\hline				
{Static-U}&	 (*)239&	-10.536	& (*)2314	&-9.0309\\ \hline
{Static-D}&188	&-11.740&		2096&	-10.536\\ \hline
{Static-DD}&	161&	-13.002&		1688&	-11.7138\\ \hline\hline
 {Heap-U}	&509&	-11.138	&	 (*)5992&	-9.7579\\ \hline
{Heap-D}	&467	&-11.740&		(*)7509 &	-10.536\\ \hline
{Heap-DD}&	450&	-12.535	&	 (*)5849&	-11.6019\\ \hline
{Heap-H}	&440	&-11.138&		 (*)5978	&-9.93399\\ \hline
{Heap-HH}	& (*)403&	-11.439	& (*)4999&	-10.235\\ \hline\hline					
{SHeap-U}&	80&	-11.740	&	 (*)717&	-11.1381\\ \hline
{SHeap-D}&	70&	-11.439	&	661	&-11.4391\\ \hline
{SHeap-DD}&	67&	-13.107&		604&	-12.3703\\ \hline
{SHeap-H}&	86	&-12.041&		 (*)662&	-10.8371\\ \hline
{SHeap-HH}&	60&	-11.689&		595&	-10.536\\ \hline		
   \end{tabular}
  \caption{{Performance of procedure {\tt SystemSolver} on the 15 models on {\tt DS1} and {\tt DS2}. The results labeled with $(*)$ are those where {\tt TFQMR} has been applied since the required precision of $10^{-11}$ on the residual norm was not satisfied after 100 steps of {\tt BCGStab}.} }\label{tab3}}

\end{center}
\end{table}

\section{Numerical Experiments}\label{numexp}

The problem of validating a ranking model is rather a difficult task since no
ground truth is known in the general case. Moreover the validity of a model
clearly depends on what we would like to measure. For example, if we want to
measure the aptitude of a scholar to work in a team we will highly value
the articles written in collaboration while if we want to measure the
scientific strength and personal skills, we may want to normalize each of the
articles by the number of co-authors.
In this respect the extreme variety of our models and the different weighting strategies allows to tune the parameters to better satisfy the different needs.

Table~\ref{tab4} summarize the experiments we performed on the two patents  datasets. In the first set of experiments we compare the different  ranking scores obtained with our models with simpler ranking methods, namely the Pagerank algorithm applied only to the citation matrix $C$, the ranking provided by one-class model and the simple citation count. The evaluation measure $P@N$ is also presented for comparing the top $N$ ranked items by some of our models with simple citation count and PageRank. The top $N$ firms obtains with some of our ranking methods are compared with the rank induced by number of patents issued by each form.

A second set of tests aims at showing that our different models are adequate to deal with incomplete data. In order to empirically prove that, we remove increasingly percentages of the attributes links to show that when dealing with incomplete database, our methods are still robust in providing a ranking ``similar'' to the one obtained with the full data. Of course, when the majority of the links are removed the rank should converge to the rank obtained with the {\tt One-class} model. A direct comparison between the top ranked results with full and partial data is done as well.

With the third set of experiments we compare the ranking scores of the same algorithms with a finer or coarser aggregation in subclasses.

 \begin{table}[t!]
\begin{center}
\small{\begin{tabular}{||l||c||r||r|| } \hline
{\tt Purpose}&{\tt experiment} & {\tt models} &{\tt Section} \\ \hline \hline
\multirow{4}*{Comparison (pat. and firms)}
&  One-class&  All& \multirow{4}*{\ref{modcomp}}\\
&  PageRank & All &  \\
& $\#$ Citations & All & \\
& P@N & {\small StiffD, StaticDD,HeapHH, SHeapHH}& \\
\hline
\multirow{3}* {Incomplete data }
&$p=0.1$&  All& \multirow{3}*{\ref{convergence}}\\
&$p=0.5$& All&  \\
&Top $N$ vs $p$ &  StaticDD&\\
\hline
\multirow{2}* {Consistence for class aggregation}
&finer& All & \multirow{2}* {\ref{class_aggr}}\\
&coarser& All&\\
\hline\hline
  \end{tabular}
  \caption{{Description of the experiments performed.}}\label{tab4}}
\end{center}
\end{table}

\subsection{Comparison between models} \label{modcomp}

The experiments reported in this section have different purposes. First we
compare the rank provided by each model with the rank obtained with the
one-class model, with the standard PageRank model and with the simple in-link
counting. The idea is that the provided rank should differ substantially from
the ranking obtained by simply counting the number of citations received, but
the presence of the features should refine the ranking without completely
reversing the importance of the players obtained by the one-class model or by
the PageRank model.

\begin{figure}[htb!]
\centering%
\includegraphics[width=7cm]{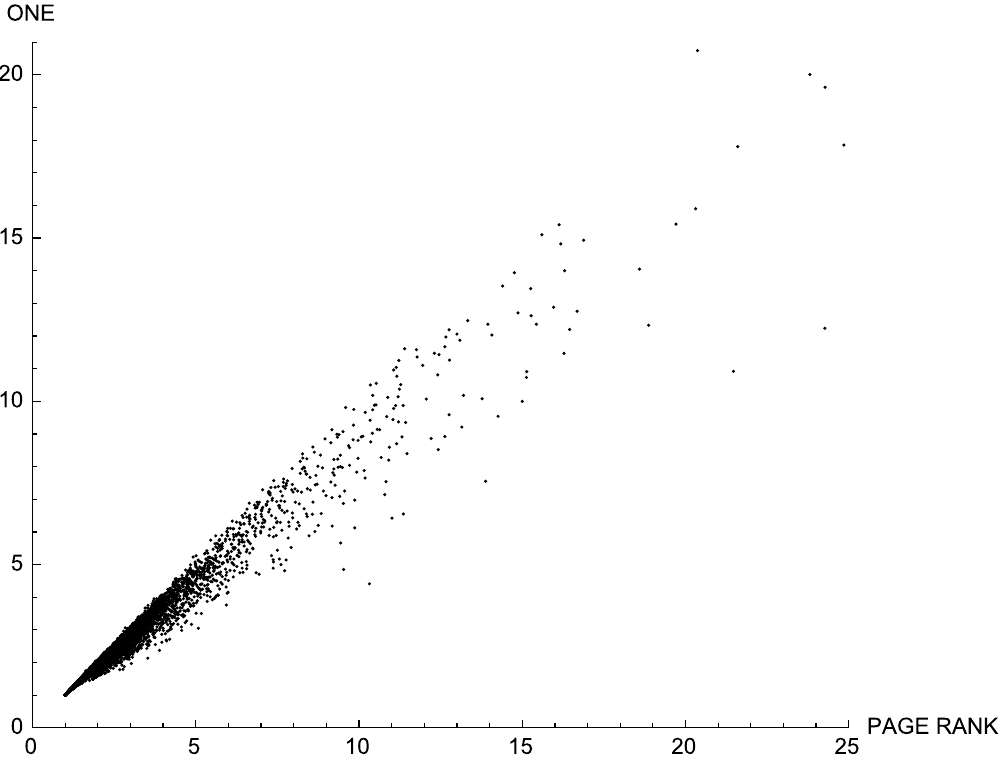}
\caption{{Comparison of the rank provided by the PageRank algorithm with
random jump probability equal to $0.15$ and the one obtained by the {\tt one-class} model
applied to DS1.}} \label{ONE-PR}
\end{figure}

In Figure~\ref{ONE-PR} it is shown the rank provided by our one-class model
versus the rank provided by the standard PageRank algorithm~\cite{BP98}.
A dot  with coordinates $(x_i, y_i)$ represent the $i$-th patent and $x_i$ is the ranking score computed with the classical PageRank algorithm, and $y_i$ the ranking score computed using our One-class model.  We
see that the two ranks are very alike because most of the points are located
on a narrow strip along the main diagonal, reflecting the high correlation between the two ranks. In fact the only difference in the
two models is the probability of reaching the dummy node which is 0.15 in the
PageRank while it changes accordingly with the outdegree of each node in our
model.

\begin{figure}[htb!]
\centering%
\includegraphics[width=14cm]{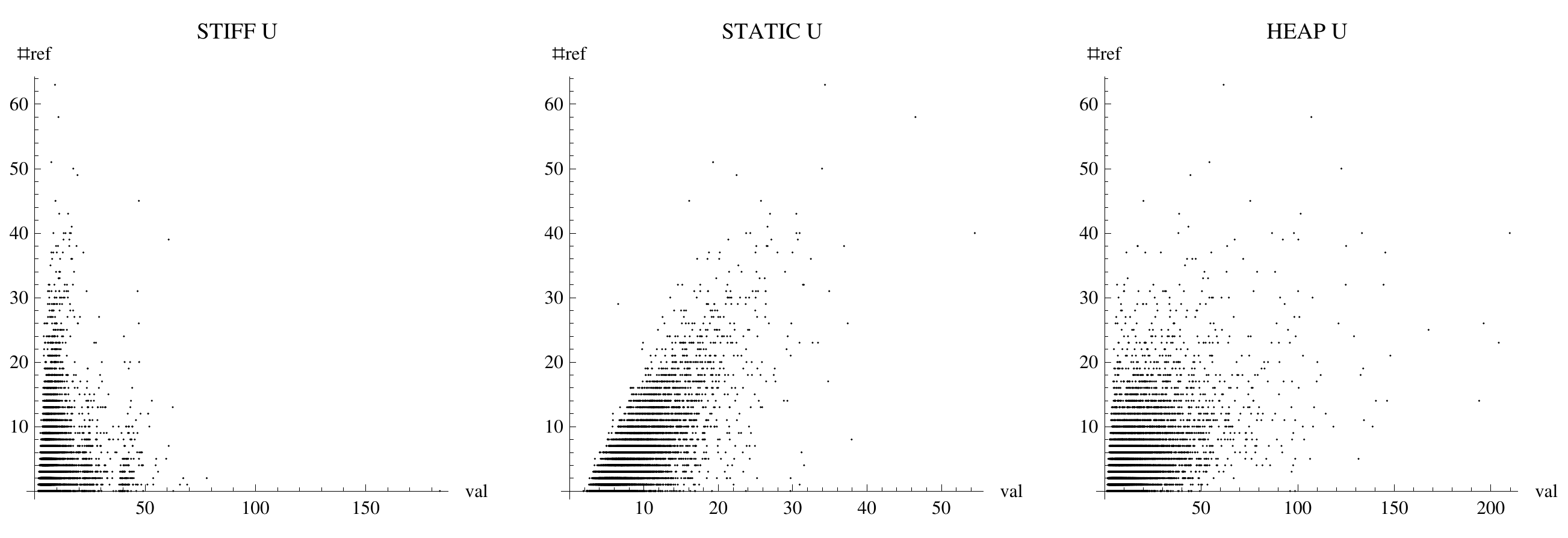}
\caption{{Comparison between the rank provided by three models with an Uniform
weighting strategy ({\tt val}) and citation count ($\#${\tt ref}).}} \label{figU}
\end{figure}
\begin{figure}[htb!]
\centering%
\includegraphics[width=14cm]{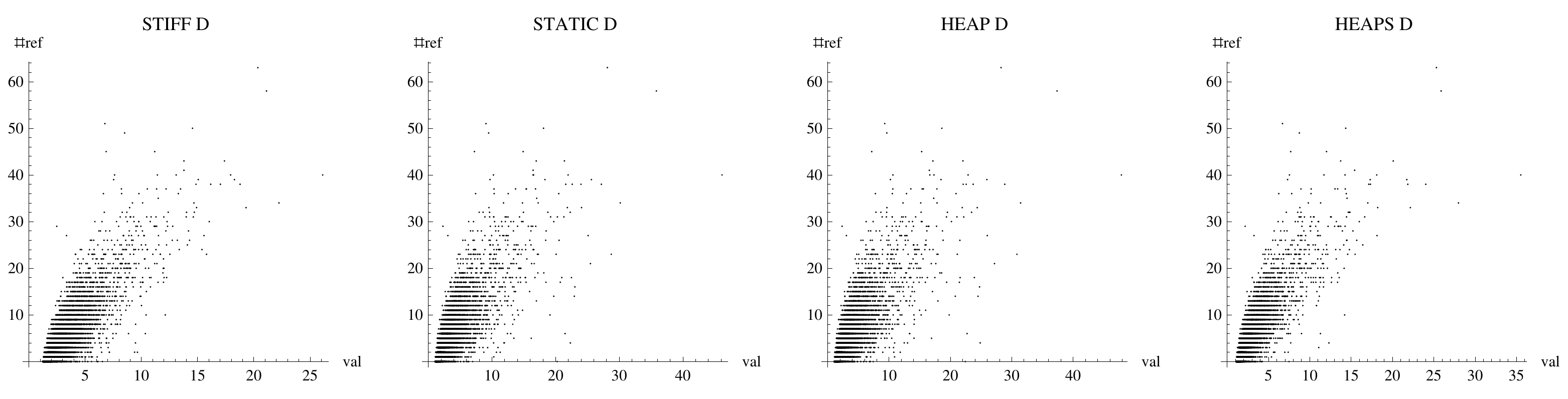}
\caption{{Comparison between the rank provided by three models with a
dimensional based weighting strategy ({\tt val}) and citation count ($\#${\tt ref}).}}
\label{FigD}
\end{figure}

Examining the plots of the ranks obtained with all the models in
Table~\ref{metodi} versus the number of citations received it turns out  that
the Uniform weighting scheme is not very adequate. In fact, for example in
Figure~\ref{figU}(a), we see that there are objects that rank very high and
have very few citations while some of those with many citations receive a
very low rank value. This effect is less noticeable in the {\tt Static} or {\tt Heap}
models but still the influence of number of citations on the actual ranking
seems to be too weak. These problems together with the instability observed
in previous section (see Tables~\ref{tab1}, \ref{tab2}, \ref{tab3} noting
that for each model procedure {\tt SystemSolver} performs better with other
weighting schemes) suggest that uniform weighting strategies are inadequate.

The results provided by most of the models using a dimension based weighting
scheme appear to be better. In fact, documents with a high number of
citations receive a good ranking score but the rank provided is not simply a
citation count. As we can observe in Figure~\ref{FigD} there is not a
substantial difference in the shape of the cloud of points obtained using
different models. Similar results can be observed with double-dimension or
heap weighting strategy.

\begin{figure}[ht]
\begin{minipage}[b]{0.45\linewidth}
\centering
\includegraphics[width=\textwidth]{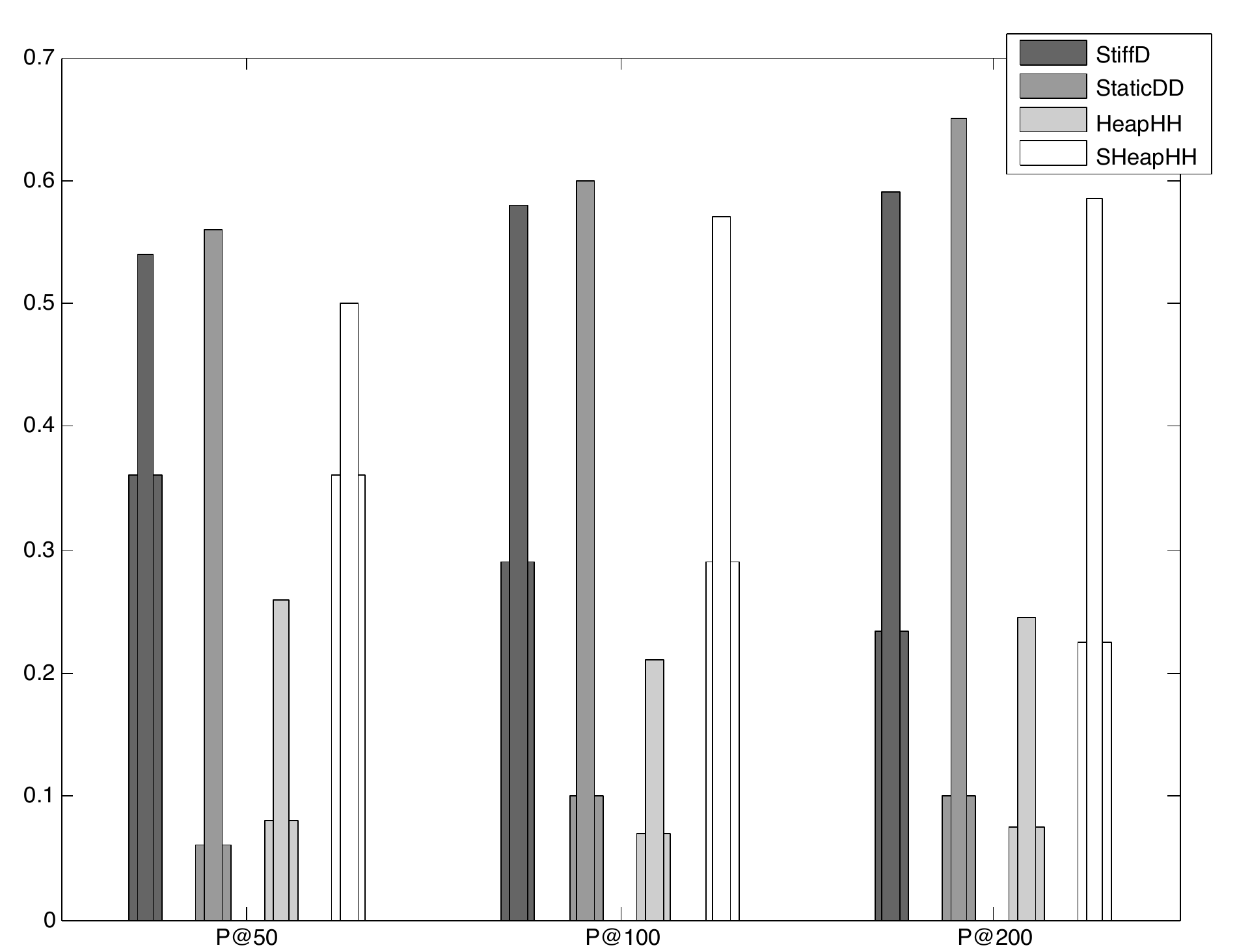}
\caption{{P@N performance of four of our models, i.e. {\tt Stif-D, Static-DD, Heap-HH, SHeap-HH} respect to citation count
(for each color the thick bars ) and PageRank (for each color the thin bars)}} \label{P@N}
\end{minipage}
\hspace{0.5cm}
\begin{minipage}[b]{0.45\linewidth}
\centering
\includegraphics[width=\textwidth]{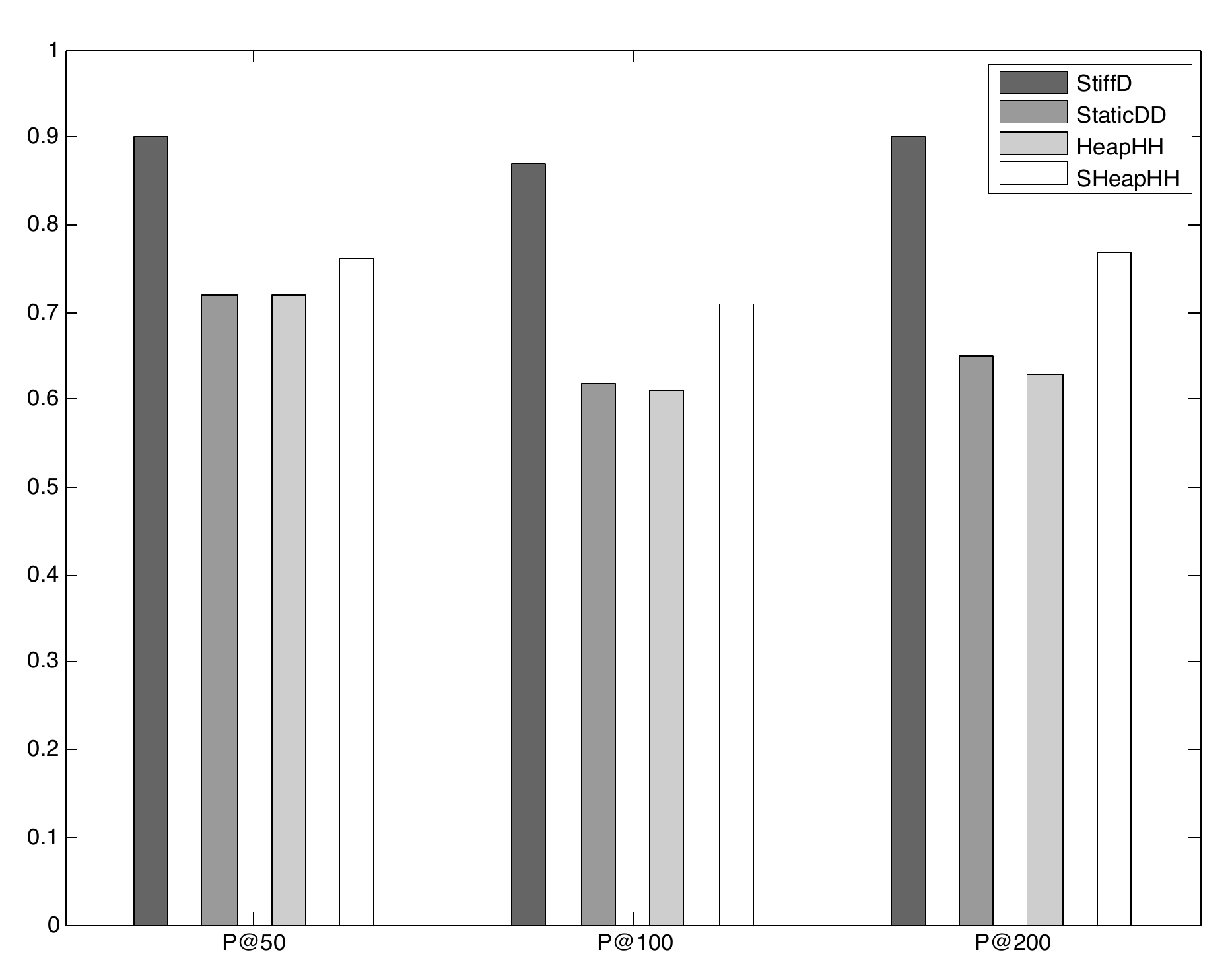}
\caption{{For firms the P@N performance of four of our models, i.e. {\tt Stiff-D, Static-DD, Heap-HH, SHeap-HH} respect to number of patent granted to each firm shows an hight correlation.}}
\label{P@Nfirm}
\end{minipage}
\end{figure}

Many authors use the  {\it precision-at-$N$} ($P@N$) measure as evaluation method. This measure is defined as follows, for a given $N\in \mathbb{N}$
$$
P@N=\frac{|E_N\cap F_N|}{N},
$$
where $E_N$ are the top ranked $N$ objects according to the ranking method one has to evaluate, and $F_N$ are the top ranked $N$  objects accordingly with the ``perfect'' ranking. Of course since the ``perfect'' ranking is not available, the top objects are generally manually ranked by volunteers or other algorithms are taken into consideration. In our case, when ranking patents it is very hard to find reliable volunteers because of the expertise required to find the most valuable patents into a such large database. We used instead as comparison the rank provided by PageRank and the citation count.
Figure~\ref{P@N} for values of $N=50, 100, 200$ depicts the performance of four of our models, i.e. {\tt Stiff-D, Static-DD, Heap-HH, SHeap-HH} respect to citation count
(thick bars) and PageRank (thin bars).
We note that our methods are more related to the rank produced by PageRank than to the simple citation count. The similarity is higher for the {\tt SHeap} model since in that case the attributes are used in a less significant way. Surprisingly enough the {\tt Static-DD model} shares more of 60\% of the top hits with PageRank, despite the two models are very different.

The precision measure $P@N$ can be used also to evaluate the firms. In Figure~\ref{P@Nfirm} we show the comparison with the rank induced by sorting the firms by the number of patents issued. We see that there is a very high correlation with the number of patents issued by a given firm, up to $90\%$ for the  {\tt  Stiff-D} model. The precision is lower for the {\tt Heap-HH} model where the citations matrix is combined with those of the features mitigating the effect of the the number of patents granted by a firm. In all the models in the top position we find very popular firms such as: IBM, Canon, Motorola, Philips, Sony, Bell {\it etc.}. Among the top results we have also firms such as Bell Labs, or Bayer AG, that despite in the time range [1976-2012] have issued a relatively low number of patents (2,617 and 896 respectively) show at the top of the list.


\subsection{Convergence with incomplete data}\label{convergence}

An important problem when dealing with large collections of multivariate data
is the incompleteness of the data. To see how robust our methods are when
part of the data are missing,  we performed many experiments leaving the
citation matrix unaltered and varying the level of information about the
features. In particular, we construct feature matrices $\tilde F_s$ obtained
taking a nonzero from $F_s$ with a fixed probability $p$, that is
$$
P\left(\tilde F_s(i,j)=1\right)=p\, F_s(i,j).
$$
Then we replace in all the models the matrices $F_s, s=1, \ldots, f$
with the matrices $\tilde F_s$.

\begin{figure}[htb!]
\centering%
\includegraphics[width=14cm]{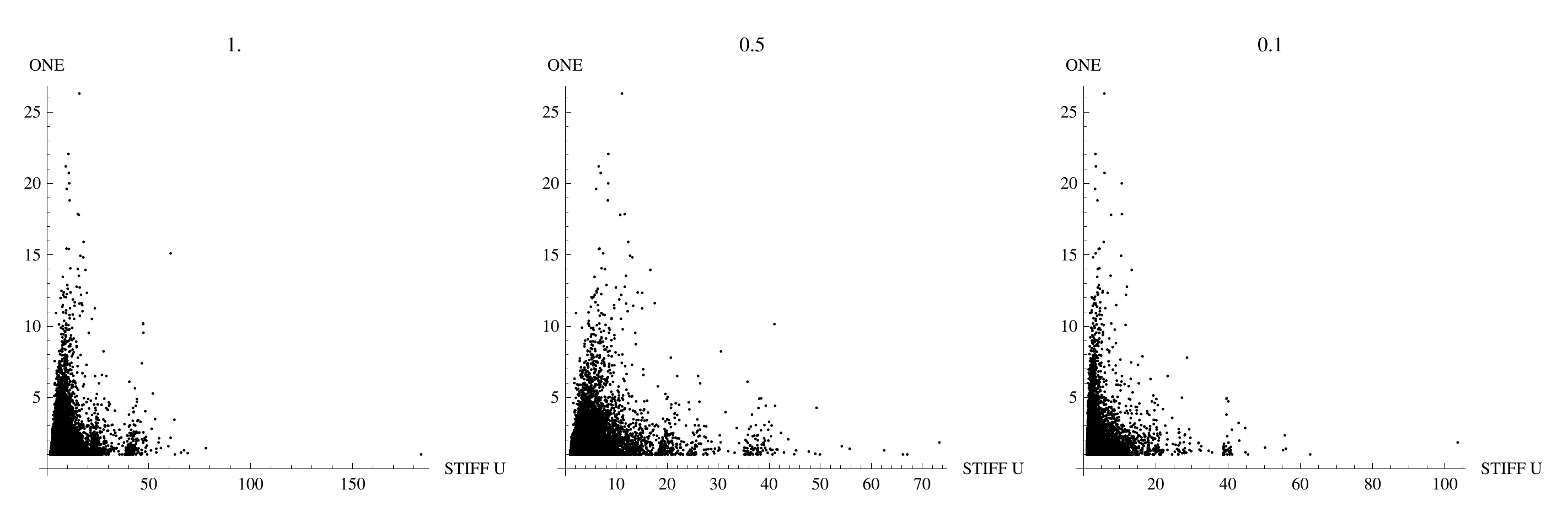}
\caption{{Comparison between the rank provided by the one-class  model and the {\tt Stiff-U} model for different values of the probability $p$.}}
\label{pstiffu}
\end{figure}

The experiments performed have two different purposes. First, we would like
to test if there are models for which the rank obtained decreasing the number
of nonzero in the feature matrices does not converge to the one obtained with
the one-class model. In fact a good model should exhibit a smooth convergence
to the one-class model as $p$ goes to zero. Second, we are interested to see
if some of the models are predictive, in the sense that the rank obtained
with missing data is ``close enough'' to the rank obtained using the full
data, suggesting a good behavior when the data are partial or missing.

\begin{figure}[htb!]
\centering%
\includegraphics[width=14cm]{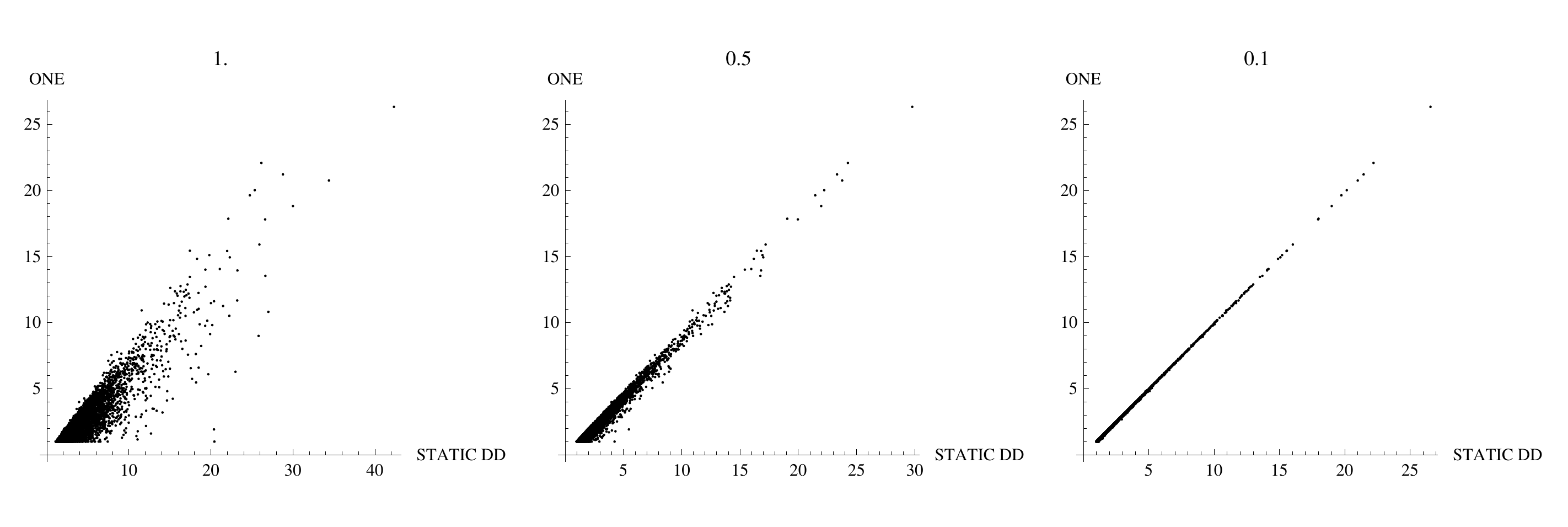}
\caption{{Comparison between the rank provided by the one-class  model and the {\tt Static-DD} model for different values of the probability $p$.}}
\label{pstaticdd}
\end{figure}

\begin{figure}[htb!]
\centering%
\includegraphics[width=14cm]{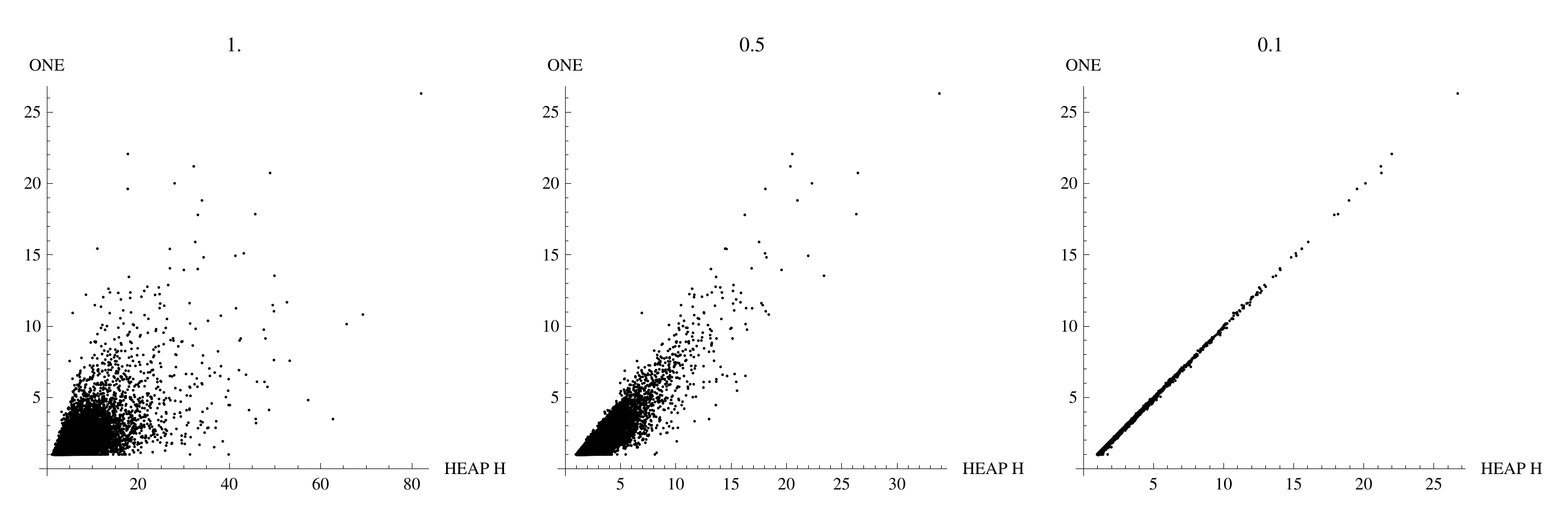}
\caption{{Comparison between the rank provided by the one-class  model and the {\tt Heap-H} model for different values of the probability $p$.}}
\label{pheaph}
\end{figure}

We report some plots obtained for values of $p$ equal to 1, 0.5 and 0.1. For
$p=0.1$ only $10\%$ of the attributes are present so the ranking obtained
should be very similar to the one obtained using only citations. Plotting the
ranking values versus the rank obtained with the one-class model, we see that
Uniform weighting schemas behave very poorly, since there is no convergence
(see Figure~\ref{pstiffu}). This fact, confirms the observation in the
previous section about the inadequateness of Uniforms weighting strategies.
On the contrary with the other weighting schemas all models exhibit a good
convergence, showing the robustness to missing data. In
Figure~\ref{pstaticdd} and ~\ref{pheaph} are depicted the results for the
three values of $p$ for the {\tt Static-DD} and the {\tt Heap-H} models.

\begin{figure}[ht]
\begin{minipage}[b]{0.45\linewidth}
\centering
\includegraphics[width=\textwidth]{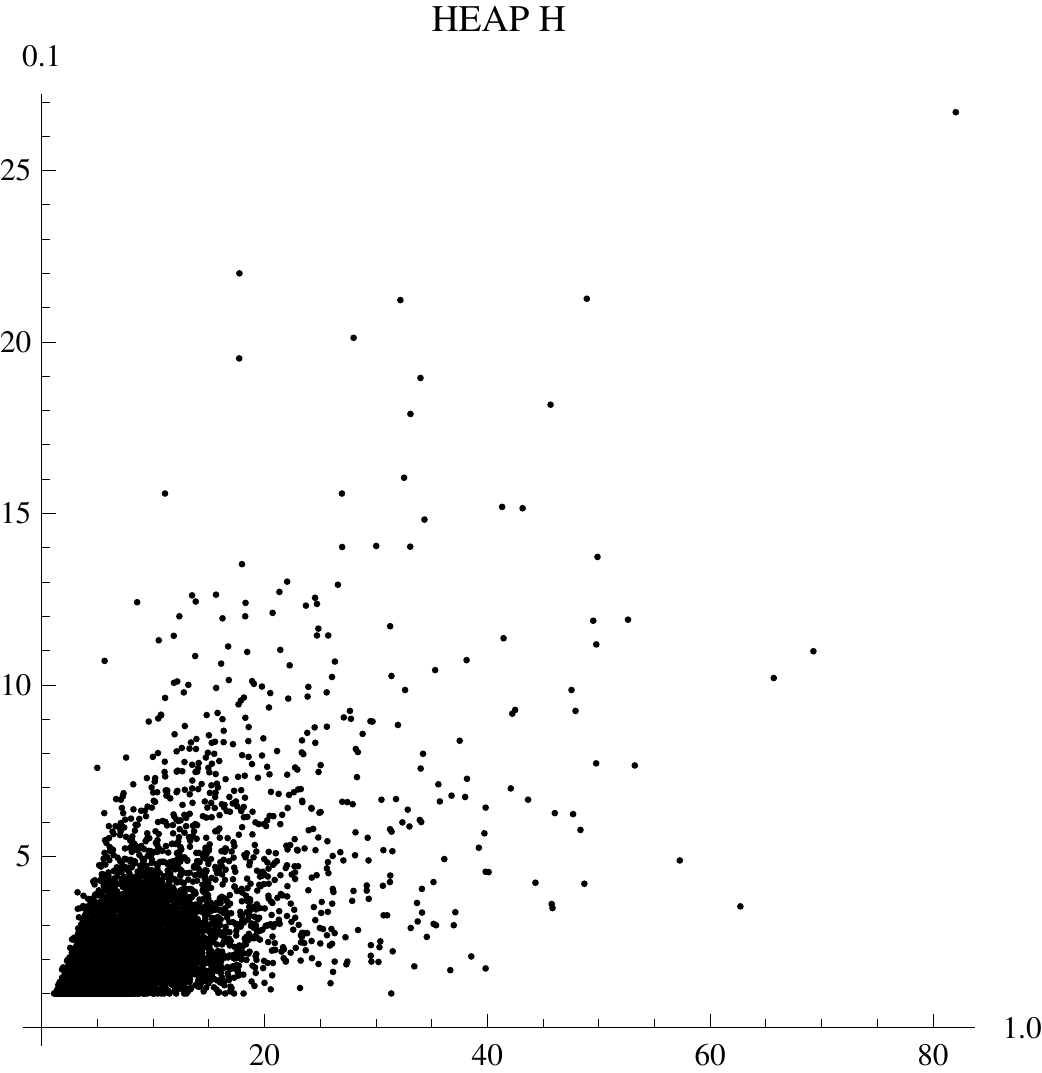}
\caption{{Comparison between the rank provided by the {\tt Heap-H} model for the patents and the same model using only 10\% of the links.}}
\label{crossheaph1}
\end{minipage}
\hspace{0.5cm}
\begin{minipage}[b]{0.45\linewidth}
\centering
\includegraphics[width=\textwidth]{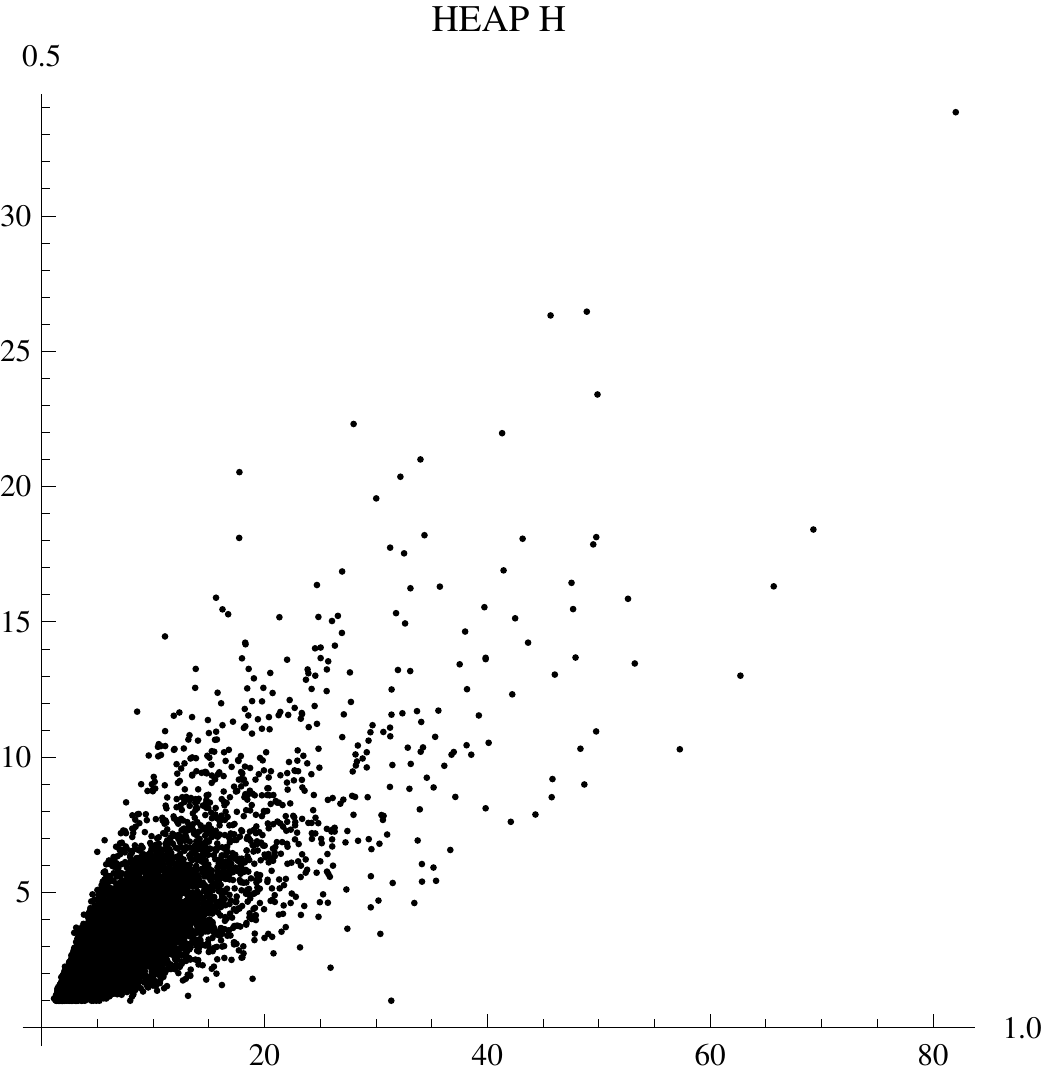}
\caption{{Comparison between the rank provided by the {\tt Heap-H} model for the patents and the same model using only 50\% of the links.}}
\label{crossheaph5}
\end{minipage}
\end{figure}

To better understand the effectiveness of the proposed methods when links are missing, we can compare the rank provided with all the links with that obtained using a small percentage of the link of the  features. In Figure~\ref{crossheaph1} and~\ref{crossheaph5} are depicted the comparison between the rank of the patents for dataset DS1, and the rank of the patents using only 10\% or 50\% of the links of the features for the {\tt Heap-H} model. We see that the rank obtained with partial information are not the same of those provided using the full matrix, but however  the cloud has a reasonable shape, showing a good predictive properties of these models for missing data. Moreover, for lower percentage of missing links, the cloud is located in an thinner region around the diagonal.

For {\tt Static-DD} model, and for the first $N$ position in the ranked list, we measure the intersection between the rank provided with the full data and the one obtained with only $10\%$ or $50\%$ of the links of the attributes. The results  in Table~\ref{tabmis} show that the rank of the patents are very similar since among the top 100 patents we have that 77 are still in the top position even removing 50\% of the links of the features, meaning that the most interesting patents show in the top position also with incomplete data.
 \begin{table}[t!]
\begin{center}
\begin{tabular}{|r||l|l|l||lll||| } \hline
{\tt N} & p=0.1&p=0.5  \\ \hline \hline
50 &  62\% & 66\% \\
100  & 74\% & 77\%\\
200 & 73\%& 77\%\\
\hline
  \end{tabular}
  \caption{{Measure of intersection between the top $N$  patents  ranked using the {\tt Static-DD} model and the rank obtained with the same model removing each edge of the attributes with probability $0.1$ or $0.5$).   }}\label{tabmis}
\end{center}
\end{table}

\subsection{Consistence for class aggregation}\label{class_aggr}

For some problems it is possible to tune the granularity of the subdivision
in classes. For example, in our databases of patents we can decide how to
group the technologies (in classes or subclasses) or geographical areas
(regions or nations) and for scientific publications we can classify papers
on the basis of their specific subject classification (there are many subject
classifications tables such as AMS, MSC,  ACM) or use a coarser grain based
on disciplines. The granularity chosen depends of course on what the ranking
is used for, but a good ranking schema should provide compatible results when
using different granularities.

\ignora{In particular we expect that summing together the rank values
obtained for entities in all the subclasses we get values similar to those
obtained working with a skinner matrix where all the subclasses are collapsed
in a unique class.}

As an example, consider the patents in Table~\ref{tec}. All these patents are
in the same  class 15 (BRUSHING, SCRUBBING, AND GENERAL CLEANING) but have a
secondary subclass as well. The rank of the patents obtained using the
extended Technologies--Patent matrix should be similar to that obtained using
a more compact Technologies--Patent matrix where, for example, the four
patents in Table~\ref{tec} are all grouped under the same Technology 15.

\begin{table}
\begin{center} { \begin{tabular}{||l||l||l||}\hline
{\tt Patent number} &{\tt Technology}& {\tt Subclass}\\ \hline \hline
6895624 & 15 &111 Brush and scraper\\ \hline
6895625&	15	&28 Rotary disk\\ \hline
6895626&	15     &50.1 Scrubber\\ \hline
6895627&	15     &98 Floor and wall cleaner:\\ \hline
   \end{tabular}
  \caption{{Four patents in the class 15- BRUSHING, SCRUBBING, AND GENERAL CLEANING, with different  subclasses.}}\label{tec}} 
\end{center}
\end{table}

\begin{figure}[ht]
\centering
\includegraphics[width=\textwidth]{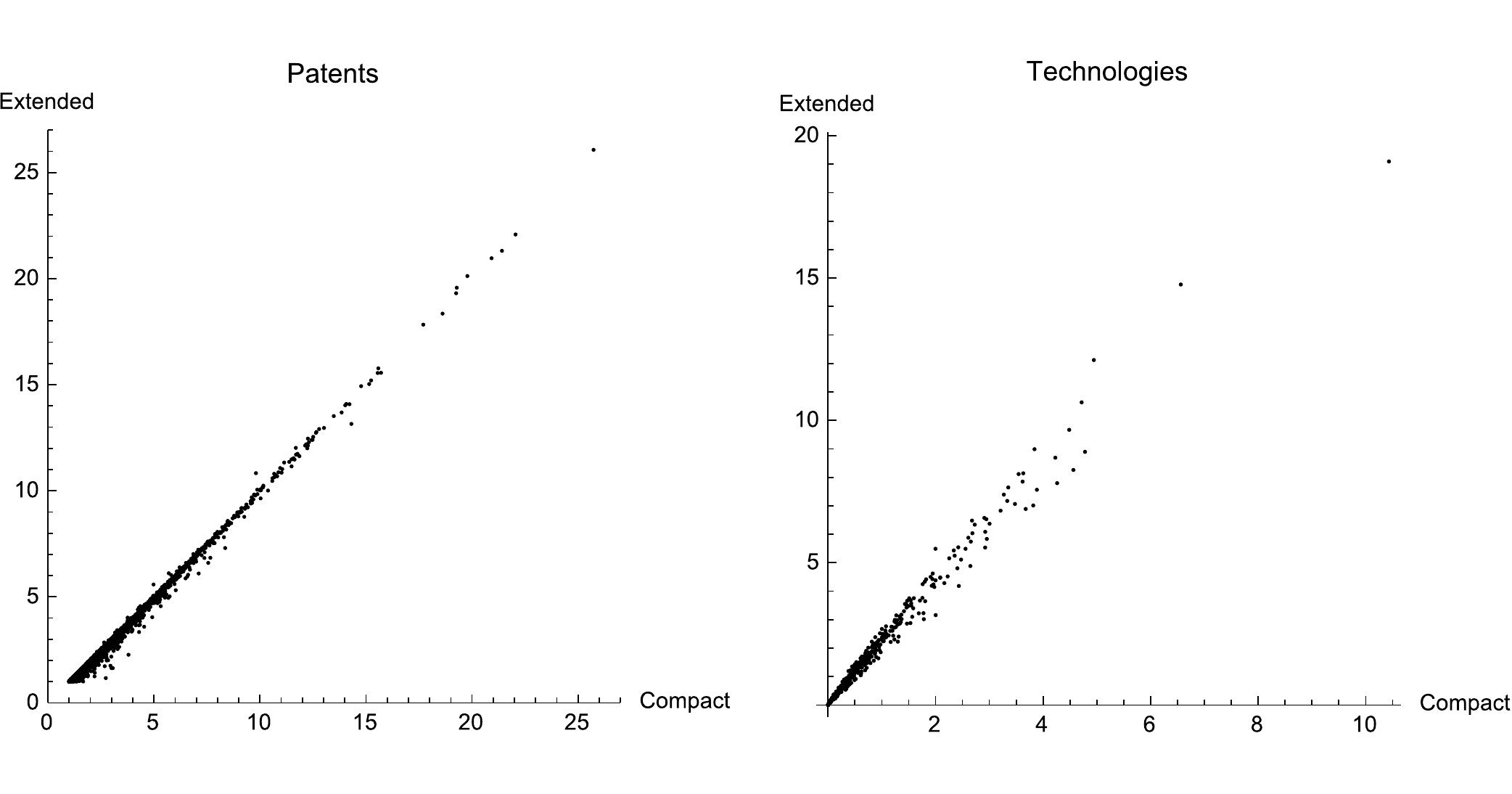}
\caption{{In the first picture a comparison between the rank of the patents provided by the {\tt Static-D} model using the extended and compact Technology--Patent matrix. In the second picture the comparison of the ranks of Technologies for the extended and compact model.}}
\label{staticDpat}
\end{figure}

Figure~\ref{staticDpat} shows the comparison of the patents' ranks obtained
using two different Technology-Patent matrices. In the compacted model we use
only the main technology class, i.e. in the example of Table~\ref{tec} the
four patents associated with different subclasses will be classified as
belonging to the same class 15. In the extended model, on the contrary, we
will use a fatter Technology-Patent matrix, with a row for each different
subclass. We see that the rank of the patents is minimally affected by the
change. Of course the rank of Technologies changes a bit more. To compare the
rank of the main 472 technologies (compact model) we summed up the rank of
all the subclasses (extended model) of a technology.

The plot obtained using models with a weighting scheme of type {\tt DD} are
less grouped around the diagonal, meaning that the ranks obtained with the
compact or extended technologies differ more than using the dimension based
technique.

\subsection{Considerations about the execution time}\label{time}

Our ranking algorithm works offline,  in the sense that the scores are precomputed and stored as is done for Google's ranking algorithm. The computation of the rank can be done periodically, for example for patents, it is reasonable to update the rank weekly, while for scientific papers a recompilation after a month would be sufficient since most of the journals have monthly issues. Our algorithms require a time ranging from 20 minutes to 2 hours to compute the ranking on the larger dataset DS2 (where the matrix involved has size approximately  of 13 millions) on a quad-core Intel Xeon @2.8GHz. To search among the documents one has to add a  search module and retrieve the documents relevant to a given query. The ranking score of the relevant documents can be simply obtained pulling out from the list of all the documents sorted by ranking score.

\section{Conclusion}\label{conclusions}

In this paper we propose several models for ranking multi-parameters data on
the basis of the linkage structure. We assume the citation matrix is enriched
with other attributes (features) that can be represented by multi-class
models. We use the attributes to improve the ranking process and as a by
product we obtain a ranking of the attributes as well. After describing the
models and different weighting strategies for measuring the influence of each
feature in the ranking process, we describe an algorithm for computing the
rank based on an iterative scheme which combines non-stationary and
stationary methods. We test some of the numerical methods on two large
datasets of US patents, and we address issues such as stability and
convergence of the algorithm applied to each model, convergence with
incomplete data, and consistence for class aggregation. In particular, the
experimental part on large datasets shows that these techniques can be used
in real applications where we have objects with multiple attributes and where
some information can be missing due to the errors or incompleteness of the
data. To search among the ordered list of objects one has to add a searching
module and retrieve only objects relevant to a given query in analogy to what
is done in the context of web search engines.


\ignora{These experiments allowed us to tune a good and effective algorithm
showing that most of the models proposed give reasonable results on the
dataset of patent data considered. In particular, we show that most of the
models have some robustness in the case of missing data by comparing results
for the complete and partial datasets. Moreover, we show that the granularity
of aggregations in feature classes does not substantially change the results,
in the sense that the resulting ranks are similar for compact or extended
aggregations in subclasses.}

As a future research we plan to address the problem of spam introducing
mechanisms for penalizing self-citations and spammers. A possible approach to
deal with cheating could consists in appropriately weighting citations and in
modifying the main diagonal of the diagonal blocks to mitigate the influence
of spammers on the final rank.

Another challenging future work is the incorporation of a preprocessing phase
aimed at recovering missing entries. Unfortunately automatic techniques such
the one proposed in~\cite{GS09} do not seem straightforwardly applicable to
our case, but maybe an attempt to recover data based on similarity of data in
specific domains such as bibliographic data can be employed. For particular
problems, such as bibliographic ranking, static indicators for journals such
as Impact Factor or Mathematical citation quotient, are available.  We plan
to investigate how this information can be used in our scheme for improving
the ranking process or as a starting point to reduce the number of
iterations.

\section*{Acknowledgments} We would like to thank Monte J. Shaffer for the many
discussions and for providing us the US patent office data. Many thanks also
to the anonymous reviewers whose comments greatly improved the quality of the
presentation.


\begin{thebibliography}{10}

\bibitem{ADKM11}
E.~Acar, D.~M. Dunlavy, T.~G. Kolda, and M.~M{\o}rup.
\newblock Scalable tensor factorizations for incomplete data.
\newblock {\em Chemometrics and Intelligent Laboratory Systems}, 106(1):41--56,
  March 2011.

\bibitem{Al09}
Paul~D. Allison.
\newblock Missing data.
\newblock In Roger~E. Millsap and Alberto Maydeu-Olivares, editors, {\em The
  SAGE: Handbook of Quantitative Methods in Psychology}, volume~3, pages
  72--89. SAGE pub, 2009.

\bibitem{Ar62}
K.~Arrow.
\newblock Economic welfare and the allocation of resources for invention.
\newblock In {\em The Rate and Direction of Inventive Activity: Economic and
  Social Factors}, Universities-National Bureau, pages 609--626. Princeton
  University Press, 1962.

\bibitem{BESS12}
G.~Berardi, A.~Esuli, F.~Sebastiani, and F.~Silvestri.
\newblock Endorsements and rebuttals in blog distillation.
\newblock {\em Information Sciences}, 249:38--47, 2013.

\bibitem{BDR}
D.~A. Bini, G.~M. {Del Corso}, and F.~Romani.
\newblock Evaluating scientific products by means of citation-based models: a
  first analysis and validation.
\newblock {\em Electron. Trans. Numer. Anal.}, 33:1--16, 2008\slash 2009.

\bibitem{BDR2}
D.~A. Bini, G.~M. {Del Corso}, and F.~Romani.
\newblock A combined approach for evaluating papers, authors and scientific
  journals.
\newblock {\em J. Comput. Appl. Math.}, 234:3104--3121, October 2010.

\bibitem{BF12}
E.~Bozzo and D.~Fasino.
\newblock A multiparameter model for link analysis of citation graphs.
\newblock {\em Electronic Transactions on Numerical Analysis}, 39:464--475,
  2012.

\bibitem{BP98}
S.~Brin and L.~Page.
\newblock The anatomy of a large-scale hypertextual {Web} search engine.
\newblock {\em Computer Networks and ISDN Systems}, 30(1--7):107--117, 1998.

\bibitem{CHT09}
F.~Chung, P.~Horn, and A.~Tsiatas.
\newblock Distributing antidote using pagerank vectors.
\newblock {\em Internet Mathematics}, 6(2):237--254, 2009.

\bibitem{Co02}
Wesley~M. Cohen, Akira Goto, Akiya Nagata, Richard~R. Nelson, and John~P.
  Walsh.
\newblock {R\&D spillovers, patents and the incentives to innovate in Japan and
  the United States}.
\newblock {\em Research Policy}, 31(8-9):1349--1367, 2002.

\bibitem{DGR07}
G.~M. {Del~Corso}, A.~Gull\'{\i}, and F.~Romani.
\newblock Comparison of {K}rylov subspace methods on the {P}age{R}ank problem.
\newblock {\em Journal of Comput. and Appl. Math.}, 210:159--166, 2007.

\bibitem{DR09}
G.~M. {Del Corso} and F.~Romani.
\newblock Versatile weighting strategies for a citation-based research
  evaluation model.
\newblock {\em Bull. Belg. Math. Soc. Simon Stevin}, 16(4):723--743, 2009.

\bibitem{DKK06}
D.~M. Dunlavy, T.~G. Kolda, and W.~P. Kegelmeyer.
\newblock Multilinear algebra for analyzing data with multiple linkages.
\newblock Technical Report SAND2006-2079, Sandia National Laboratories,
  Albuquerque, NM and Livermore, CA, April 2006.

\bibitem{GZB04}
D.~F. Gleich, L.~Zhukov, and P.~Berkhin.
\newblock Fast parallel {PageRank}: A linear system approach.
\newblock Technical Report YRL-2004-038, Yahoo! Research Labs, 2004.

\bibitem{GS09}
Roger Guimerˆ and Marta Sales-Pardo.
\newblock Missing and spurious interactions and the reconstruction of complex
  networks.
\newblock {\em Proceedings of the National Academy of Sciences},
  106(52):22073--22078, 2009.

\bibitem{Kl99}
J.~M. Kleinberg.
\newblock Authoritative sources in a hyperlinked environment.
\newblock {\em Journal of the ACM}, 46(5):604--632, 1999.

\bibitem{KB06}
T.~Kolda and B.~Bader.
\newblock The {TOPHITS} model for higher-order web link analysis.
\newblock In {\em Proceedings of Link Analysis, Counterterrorism and Security
  2006}, 2006.

\bibitem{LR87}
R.~J.~A. Little and D.~B. Rubin.
\newblock {\em Statistical Analysis with Missing Data}.
\newblock John Wiley \& Sons, New York, NY, USA, 1987.

\bibitem{NZW05}
Z.~Nie, Y.~Zhang, J.~Wen, and W.~Ma.
\newblock Object-level ranking: bringing order to web objects.
\newblock In {\em Proceedings of the 14th international conference on World
  Wide Web}, pages 567--574, New York, NY, USA, 2005.

\bibitem{rev1}
Therese~D. Pigott.
\newblock A review of methods for missing data.
\newblock {\em Educational research and Evaluation}, 7(4):353--383, 2001.

\bibitem{Sa03}
Y.~Saad.
\newblock {\em Iterative methods for sparse linear systems (2nd Edition)}.
\newblock Society for Industrial and Applied Mathematics, Philadelphia, PA,
  USA, 2003.

\bibitem{Sh97}
J.~L. Schafer.
\newblock {\em Analysis of Incomplete Multivariate Data}.
\newblock Chapman \& Hall, London, 1997.

\bibitem{Sc02}
Judi Scheffer.
\newblock Dealing with missing data.
\newblock {\em Research Letters in the Information and Mathematical Sciences},
  3:153--160, 2002.

\bibitem{Monte1}
M.~Shaffer.
\newblock {\em Entrepreneurial Innovation: {P}atent {R}ank and {M}arketing
  science}.
\newblock PhD thesis, Washington State University, USA, 2011.

\bibitem{Su08}
J.~Suchal.
\newblock Enhancing search using layered graph ranking of multigraphs.
\newblock In {\em IIT.SRC}, pages 1--8, 2008.

\bibitem{SuHa12}
Y.~Sun and Han J.
\newblock Mining heterogeneous information networks: A structural analysis
  approach.
\newblock {\em SIGKDD Explorations}, 6(2):20--28, 2012.

\bibitem{SuYuHa09}
Yizhou Sun, Yintao Yu, and Jiawei Han.
\newblock Ranking-based clustering of heterogeneous information networks with
  star network schema.
\newblock In {\em Proceedings of the International Conference on Knowledge
  Discovery and Data Mining}, pages 797--806, 2009.

\bibitem{Va00}
R.~S. Varga.
\newblock {\em Matrix iterative analysis}.
\newblock Springer series in computational mathematics. Springer Verlag,
  Berlin, Heidelberg, Paris, 2000.

\bibitem{ZaFe11}
Ming Zhang, Sheng Feng, Jian Tang, Bolanle Ojokoh, and Guojun Liu.
\newblock Co-ranking multiple entities in a heterogeneous network: Integrating
  temporal factor and usersÕ bookmarks.
\newblock In {\em Digital Libraries: For Cultural Heritage, Knowledge
  Dissemination, and Future Creation}, volume 7008 of {\em Lecture Notes in
  Computer Science}, pages 202--211. Springer Berlin Heidelberg, 2011.

\bibitem{ZLZ09}
T.~Zhou, Linyuan L\"{u}, and Yi-Cheng Zhang.
\newblock Predicting missing links via local information.
\newblock {\em The European Physical Journal B}, 71(4):623--630, 2009.

\bibitem{ZLZ11}
T.~Zhou, Linyuan L\"{u}, and Yi-Cheng Zhang.
\newblock Leaders in social networks, the delicious case.
\newblock {\em PLoS ONE}, 6:623--630, 2011.

\bibitem{ZZJZX}
Xiaofeng Zhu, S.~Zhang, Zhi Jin, Zili Zhang, and Zhuoming Xu.
\newblock Missing value estimation for mixed-attribute data sets.
\newblock {\em IEEE Transactions on Knowledge and Data Engineering},
  23(1):110--121, 2011.

\end{thebibliography}

\end{document}